\pdfoutput=1
\documentclass[reqno]{amsart}
\usepackage{amssymb}
\usepackage{amsmath,amssymb,amsfonts,amsthm,latexsym,amscd,epsfig,color}
\usepackage{mathrsfs}

\usepackage{a4wide}
\usepackage[T1]{fontenc}
\usepackage[latin1]{inputenc}
\usepackage{tikz}
\usepackage{appendix}
\usetikzlibrary{shapes}

\renewcommand{\(}{\left(}
\renewcommand{\)}{\right)}
\newtheorem{theo}{Theorem}
\newtheorem{prop}{Proposition}
\newtheorem{lemma}{Lemma}
\newtheorem{cor}{Corollary}
\theoremstyle{definition}
\newtheorem{df}{Definition}
\newtheorem{ex}{Example}

\theoremstyle{remark}

\newtheorem{rem}{Remark\!\!}

\newtheorem{nrem}{Remark}

\newcommand{\Ord}[1]{{\mathcal O}\left(#1\right)}
\newcommand{\rcp}[1]{\frac{1}{#1}}

\newcommand{\bt}{\begin{theo}}
\newcommand{\et}{\end{theo}}
\newcommand{\bl}{\begin{lemma}}
\newcommand{\el}{\end{lemma}}
\newcommand{\bp}{\begin{prop}}
\newcommand{\ep}{\end{prop}}
\newcommand{\bdf}{\begin{df}}
\newcommand{\edf}{\end{df}}
\newcommand{\brem}{\begin{rem}}
\newcommand{\erem}{\end{rem}}
\newcommand{\bnrem}{\begin{nrem}}
\newcommand{\enrem}{\end{nrem}}
\newcommand{\bex}{\begin{ex}}
\newcommand{\eex}{\end{ex}}
\newcommand{\bcor}{\begin{cor}}
\newcommand{\ecor}{\end{cor}}
\newcommand{\bncor}{\begin{ncor}}
\newcommand{\encor}{\end{ncor}}
\newcommand{\bpf}{\begin{proof}}
\newcommand{\epf}{\end{proof}}
\newcommand{\eps}{\varepsilon}

\newcommand{\nti}{{n\to\infty}}
\newcommand{\M}{\mathcal M}
\newcommand{\Y}{\mathbf Y}

\newcommand{\El}[1]{\small $\ell_{#1}$}
\newcommand{\Ke}[1]{$k_{#1}$}
\newcommand{\Rl}[1]{$r_{#1}$}

\begin{document}

\title[Phylogenetic networks with few reticulation vertices]{Counting phylogenetic networks with few reticulation vertices: tree-child and normal networks}
\author{Michael Fuchs \and Bernhard Gittenberger \and Marefatollah Mansouri}
\thanks{This research has been supported by a bilateral Austrian-Taiwanese project FWF-MOST, grants
I~2309-N35 (FWF) and MOST-104-2923- M-009-006-MY3 (MOST)}

\address{MF: Department of Applied Mathematics, National Chiao Tung University,
Hsinchu, 300, Taiwan.}
\email{mfuchs@math.nctu.edu.tw}

\address{BG and MM: Department of Discrete Mathematics and Geometry, Technische
Universit\"at Wien, Wiedner Hauptstra\ss e 8-10/104, A-1040 Wien, Austria.}
\email{\{gittenberger,marefatollah.mansouri\}@dmg.tuwien.ac.at}




\date{\today}

\begin{abstract}
In recent decades, phylogenetic networks have become a standard tool in modeling evolutionary processes. Nevertheless, basic combinatorial questions about them are still largely open. For instance, even the asymptotic counting problem for the class of phylogenetic networks and subclasses is unsolved. In this paper, we propose a method based on generating functions to count networks with few reticulation vertices for two subclasses which are important in applications: tree-child networks and normal networks. In particular, our method can be used to completely solve the asymptotic counting problem for these network classes when the number of reticulation vertices remains fixed while the network size tends to infinity.
\end{abstract}

\maketitle

\section{Introduction and Results}

This paper will be concerned with the counting of {\it phylogenetic networks}, a basic and
fundamental question which is of interest in mathematical biology; see \cite{MSW15}.

Phylogenetic networks are used to model reticulation events in evolutionary biology. Even though
the presence of such events has been acknowledged by biologists since the dawn of the development
of evolution as a scientific discipline, for the most part, {\it phylogenetic trees} instead of phylogenetic networks
have been used to model the relationship between species. This might be due to the fact that trees
are a considerably simpler structure than networks and thus allow a rich theory. For instance,
their combinatorics is well-understood: the corresponding counting problem was already solved by
Schr\"{o}der in 1870 \cite{Sc}. Several further studies were published to analyze parameters or
variations, e.g. \cite{Bo16, BoFl09, FoRo80, FoRo88}. Moreover, phylogenetic trees are also important for constructing
phylogenetic networks (see \cite{BoSe16a, ChWa12}) and thus the comparison of phylogenetic
trees and networks is an active area of research, see \cite{FrSt15} and \cite{BoSe16b, CLS14, Se16} concerning
tree-embeddings in networks.

The combinatorics of phylogenetic networks, on the other hand, remains a challenge and only few
papers have addressed it. It is the goal of this paper to make some further progress and in
particular to solve the counting problem for phylogenetic networks with a fixed number of
reticulation vertices.

Before stating our results in more detail, we recall some definitions and previous work. First,
a phylogenetic network is defined as a rooted directed acyclic graph (DAG) which
is connected and consists of the following vertices:  \begin{itemize} \item[(i)] a {\it root vertex}
which has in-degree $0$ and out-degree $2$ (except if the network consists of only one vertex);
\item[(ii)] {\it tree vertices} which have in-degree
$1$ and out-degree $2$; \item[(iii)] {\it reticulation vertices} which have in-degree $2$ and
out-degree $1$; \item[(iv)] {\it leaves} which have in-degree $1$ and out-degree $0$.
\end{itemize} Phylogenetic networks are usually labeled, where all labels are assumed to be
different and two kinds of labelings are often considered: (i) all vertices are labeled; such
networks we will call {\it vertex-labeled networks} throughout this work, and (ii) only leaves
are labeled; these are called {\it leaf-labeled networks}.

Phylogenetic networks are used to model reticulate evolution. However, the process of evolution
is driven by specific principles which add further restrictions on phylogenetic networks. Thus,
biologists have defined many subclasses of the class of phylogenetic networks. Two of them are
{\it tree-child networks} and {\it normal networks}; see e.g. \cite{CaRoVa,Wi}.

\begin{figure}[h]
\setlength{\unitlength}{1cm}
\begin{center}
        \begin{tikzpicture}[scale=1.1]
        \draw (0cm,0cm) node[circle,inner sep=1,fill,draw] (1) {};
        \draw (-1cm,-1cm) node[circle,inner sep=1,fill,draw] (2) {};
        \draw (0cm,-1cm) node[circle,inner sep=1,fill,draw] (3) {};
        \draw (1cm,-1cm) node[circle,inner sep=1,fill,draw] (4) {};
        \draw (-0.5cm,-1.5cm) node[circle,inner sep=1,fill,draw] (5) {};
        \draw (0.5cm,-1.5cm) node[circle,inner sep=1,fill,draw] (6) {};
        \draw (-2.1cm,-2.1cm) node[circle,inner sep=1,fill,draw] (7) {};
        \draw (-.8cm,-2.1cm) node[circle,inner sep=1,fill,draw] (8) {};
        \draw (-.2cm,-2.1cm) node[circle,inner sep=1,fill,draw] (13) {};
        \draw (0.5cm,-2.1cm) node[circle,inner sep=1,fill,draw] (9) {};
        \draw (2.1cm,-2.1cm) node[circle,inner sep=1,fill,draw] (10) {};
        \draw (-0.5cm,-0.5cm) node[circle,inner sep=1,fill,draw] (11) {};
        \draw (0.5cm,-0.5cm) node[circle,inner sep=1,fill,draw] (12) {};
        \draw (-1.5cm,-1.5cm) node[circle,inner sep=1,fill,draw] (14) {};
        \draw (-1.2cm,-2.1cm) node[circle,inner sep=1,fill,draw] (15) {};
        \draw (1.5cm,-1.5cm) node[circle,inner sep=1,fill,draw] (16) {};
        \draw (1cm,-2.1cm) node[circle,inner sep=1,fill,draw] (17) {};
        \put(0, -2.7){\makebox(0, 0){$(i)$}}
        \draw (1) -- (4); \draw (4) -- (10); \draw (1) -- (11);
        \draw (11) -- (2); \draw (11) -- (3); \draw (2) -- (5);\draw (16) -- (17);
        \draw (3) -- (12); \draw (3) -- (6); \draw (4) -- (6);\draw (14) -- (15);
        \draw (2) -- (7); \draw (5) -- (8); \draw (6) -- (9);\draw (13) -- (5);

        \draw (5cm,0cm) node[circle,inner sep=1,fill,draw] (1) {};
        \draw (4.7cm,-.3cm) node[circle,inner sep=1,fill,draw] (2) {};
        \draw (5cm,-.7cm) node[circle,inner sep=1,fill,draw] (3) {};
        \draw (4cm,-1.2cm) node[circle,inner sep=1,fill,draw] (4) {};
        \draw (6.3cm,-1.2cm) node[circle,inner sep=1,fill,draw] (5) {};
        \draw (5cm,-1.2cm) node[circle,inner sep=1,fill,draw] (6) {};
        \draw (4.5cm,-1.75cm) node[circle,inner sep=1,fill,draw] (7) {};
        \draw (3.5cm,-2.1cm) node[circle,inner sep=1,fill,draw] (8) {};
        \draw (4.5cm,-2.1cm) node[circle,inner sep=1,fill,draw] (9) {};
        \draw (5.5cm,-2.1cm) node[circle,inner sep=1,fill,draw] (10) {};
        \draw (7cm,-2.1cm) node[circle,inner sep=1,fill,draw] (13) {};
        \draw (6.3cm,-2.1cm) node[circle,inner sep=1,fill,draw] (14) {};
        \draw (5.5cm,-.9cm) node[circle,inner sep=1,fill,draw] (15) {};
                \put(5.9, -2.7){\makebox(0, 0){$(ii)$}}
        \draw (1) -- (2); \draw (1) -- (5); \draw (2) -- (3);
        \draw (3) -- (5); \draw (5) -- (13); \draw (2) -- (4);
        \draw (3) -- (6); \draw (10) -- (6); \draw (4) -- (7);
        \draw (6) -- (7); \draw (7) -- (9); \draw (4) -- (8);   \draw (14) -- (15);

        \end{tikzpicture}
\end{center}
\vspace*{0.25cm}
\caption{Two phylogenetic networks, where $(i)$ is a general network that is not a tree-child
network and $(ii) $ is a tree-child network that is not a normal network. Edges are directed downwards. \label{TN}}
\label{TC_normal}
\end{figure}
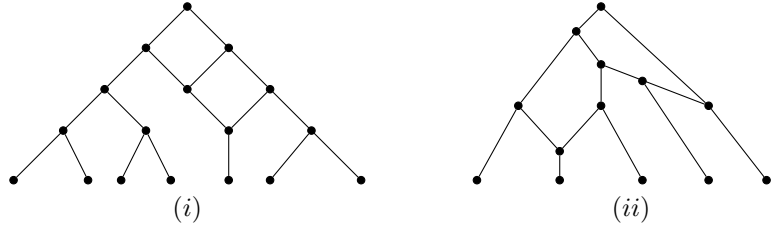

In tree-child networks, one has the additional requirement that reticulation events
cannot happen in close proximity, or more formally, every tree vertex must have at least one child
which is not a reticulation vertex and no reticulation vertex is  directly followed by another
reticulation vertex. Normal networks, on the other hand, form a subclass of the class of tree-child
networks with the additional requirement that evolution does not take shortcuts, or again more
formally, if there is a path of at least length $2$ from a vertex $u$ to a vertex $v$, then there
is no direct edge from $u$ to $v$. For examples of such networks see Figure~\ref{TC_normal}.

\brem
Note that in general phylogenetic networks, as defined above, multiple edges
are not explicitly forbidden (except when dealing with enumeration of leaf-labeled networks, since
otherwise the counting problem is not meaningful). In fact, only double edges may occur because of
the degree constraints. The tree-child condition, however, makes double edges impossible. Thus
tree-child and normal networks do not contain double edges.
\erem

Next, let us recall what is known about the number of tree-child and normal networks. Denote by $T_n$ and $N_n$
the number of vertex-labeled tree-child networks and vertex-labeled normal networks, respectively, where $n$
is the total number of vertices. Similarly, denote by $\tilde{T}_{\ell}$ and $\tilde{N}_{\ell}$ the
number of leaf-labeled tree-child networks and leaf-labeled normal networks, where $\ell$ denotes
the number of leaves. Then, it was proved in \cite{MSW15} that for all odd $n$,
\[
(e_1n)^{5n/4}\leq N_n\leq T_n\leq (e_2n)^{5n/4},
\]
where $e_1,e_2>0$ are suitable constants.
(It is easy to see that $N_n=T_n=0$ when $n$ is even.) Similarly, there are $f_1,f_2>0$ such that for all $\ell$,
\[
(f_1\ell)^{2\ell}\leq\tilde{N}_{\ell}\leq\tilde{T}_{\ell}\leq(f_2\ell)^{2\ell}.
\]
Note that the first result can be equivalently stated as
\[
N_n=n^{5n/4+{\mathcal O}(n/\log n)}\qquad\text{and}\qquad T_n=n^{5n/4+{\mathcal O}(n/\log n)}
\]
and the second as
\[
\tilde{N}_{\ell}=\ell^{2\ell+{\mathcal O}(\ell/\log\ell)}\qquad\text{and}\qquad
\tilde{T}_{\ell}=\ell^{2\ell+{\mathcal O}(\ell/\log\ell)}.
\]
Thus, one is still quite far away from getting precise asymptotics for these counting sequences
and this was left as an open problem in \cite{MSW15}.

In this paper, we will consider tree-child  and normal networks with a fixed number $k$
of reticulation vertices. It should be mentioned that they form (very) small subclasses of the
class of all tree-child and normal networks since it was also proved in \cite{MSW15}
that almost all vertex-labeled tree-child resp. normal networks have $k\sim n/4$ and almost all
leaf-labeled tree-child resp. normal networks have $k\sim\ell$. Nevertheless, these subclasses
are interesting from a combinatorial point of view since we can get precise asymptotics of their
numbers. Moreover, they are more suitable for modelling phylogenesis in environments where
reticulation is a very rare event (although even then it may be sometimes desirable to admit
$k\rightarrow\infty$ as $n\rightarrow\infty$). Models with a fixed number of reticulation vertices
were for instance considered in \cite{BoLiSe17,SeSt06}.
Likewise, in the construction of phylogenetic networks from trees, models with bounded
reticulation number do play a role, see \cite{KeSc14,vKS16}.

Recently, people studying phylogenetic networks or related structures have become more and more
interested in enumerative aspects. We mentioned already the shape analysis of phylogenetic trees
\cite{Bo16, BoFl09, FoRo80, FoRo88} and the bounds for the counting sequences of some classes of
phylogenetic networks \cite{MSW15}. But other counting problems were studied in
\cite{AlAl16,CEJM13,DiRo17,LSS13,Ro07,Se17,SeSt06}. Though combinatorial counting problems are
often amenable to the rich tool box of analytic combinatorics \cite{AnaCombi}, generating
functions have been rarely used in phylogenetic enumeration problems.

Here we focus on the already mentioned class of phylogenetic networks with a low number of
reticulation events, more specifically on the above two subclasses of this class, and demonstrate how
analytic combinatorics can be used to obtain general (asymptotic) enumeration results for those
classes. We believe that our paper is of interest to experts working on the mathematics of
phylogenetics and that many more enumeration problems in phylogenetics can be approached in a
similar way.

Now, denote by $N_{k,n}$ resp. $T_{k,n}$ the number of normal resp. tree-child networks with
$k$ reticulation vertices in the vertex-labeled case and $\tilde{N}_{k,\ell}$
resp. $\tilde{T}_{k,\ell}$ in the leaf-labeled case. Then, our results are as follows.
\bt\label{main-thm-1}
For the number $N_{k,n}$ of vertex-labeled normal networks with $k\geq 1$ reticulation vertices, there is a positive constant
$c_k$ such that
\[
N_{k,n}\sim c_k\left(1-(-1)^n\right)\left(\frac{\sqrt{2}}{e}\right)^n n^{n+2k-1},\qquad (n\rightarrow\infty).
\]
In particular,
\begin{align*}
c_1=\frac{\sqrt{2}}{4};\qquad
c_2=\frac{\sqrt{2}}{32};\qquad
c_3=\frac{\sqrt{2}}{384}.
\end{align*}
\et
\begin{rem}
Note that this result also holds for $k=0$ where it becomes the result of Schr\"{o}der; see above and \cite{Sc}.
\end{rem}
Surprisingly, the same result also holds for vertex-labeled tree-child networks. (It was proved in
\cite{MSW15} that $N_n=o(T_n)$.) This shows in particular that if one considers only first-order
asymptotics, then the additional requirement for normal networks does not matter. Note, however,
that we are considering networks with an \emph{a priori} fixed number $k$ of reticulation vertices.
Thus, we do not claim that the asymptotic equivalence given in Theorem~\ref{main-thm-1} holds uniformly in $k$ (and neither do we claim this in Theorem~\ref{main-thm-2} below). Indeed, such a claim would surely be wrong since otherwise one could sum up both sides over $k$
and would get a contradiction to the above mentioned result from \cite{MSW15}.
\bt\label{main-thm-2}
For the number $T_{k,n}$ of vertex-labeled tree-child networks with $k\geq 1$ reticulation vertices,
\[
T_{k,n}\sim c_k\left(1-(-1)^n\right)\left(\frac{\sqrt{2}}{e}\right)^nn^{n+2k-1},\qquad (n\rightarrow\infty)
\]
with $c_k$ as in the previous theorem.
\et
\begin{rem}
Again the result also holds for $k=0$ where it is Schr\"{o}der's result. Moreover, note that for
$k=0$ and $k=1$, $T_{k,n}$ is identical with the number of all vertex-labeled phylogenetic networks (for
$k\geq 2$, the latter number becomes strictly larger than $T_{k,n}$, however, the leading term of the asymptotic expansion is likely to be again the same; see Section~\ref{con}).
\end{rem}
\begin{cor}
Let $k\geq 1$. Then, asymptotically almost all tree-child networks with $k$ reticulation vertices are normal networks.
\end{cor}
\begin{rem}
When going beyond first-order asymptotics, one sees that the additional requirement for normal networks does indeed matter; see below for longer asymptotic expansions for $k=1,2,3$ which show a difference in the second order term for vertex-labeled normal and tree-child networks.
\end{rem}

Similar results to the results above will be shown for leaf-labeled tree-child and normal
networks, too; see Section \ref{ll-net}.

The remainder of  the paper is as follows. In the next section, we will explain how to use
generating functions to count tree-child and normal networks. The methodology we use is what is
known nowadays as ``Analytic Combinatorics'' \cite{AnaCombi} and relies on the symbolic
method~\cite[Sec.~I.1--I.2]{AnaCombi} and the treatment of labeled
structures~\cite[Sec.~II.1--II.2]{AnaCombi} as well as the pointing
operation~\cite[Sec.~II.6]{AnaCombi}. This counting
procedure will then be applied in Section \ref{vl-normal-net} to vertex-labeled normal networks.
(This section will contain the proof of Theorem \ref{main-thm-1}.) In Section \ref{vl-tc-net},
we apply the same approach to vertex-labeled tree-child networks and prove Theorem
\ref{main-thm-2}. In Section \ref{ll-net}, we will briefly discuss results for leaf-labeled
networks which are obtained from those for vertex-labeled networks in Section \ref{vl-normal-net}
and Section \ref{vl-tc-net}. Finally, we will conclude the paper with some remarks in Section
\ref{con}.

\section{Decomposing Phylogenetic Networks}

In order to count the above classes of phylogenetic networks, we will decompose them and use this decomposition to
obtain a reduction which can be easily analyzed by means of generating functions. Then the
reduction is extended to get back the original network in such a way that the extension procedure
has a counterpart in generating function algebra, hence allowing an asymptotic analysis of the
number of phylogenetic networks.
We start with normal networks, since tree-child networks differ from normal ones just by dropping a
condition which allows a similar analysis.

Consider a normal network having exactly $k$ reticulation vertices. Then each such vertex has
two incoming edges. Choose one of them and remove it. The remaining graph is a (labeled and
nonplane) Motzkin tree\footnote{We mention that we slightly abuse the word here: A Motzkin tree
(also known as unary-binary tree) is usually an unlabeled and plane. The concept stems from
computer science, see \cite{FlOd82,FlPr86,FlSt80}. In contrast, the trees we are considering here are
labeled and nonplane, but nevertheless still unary-binary trees. Thus, they are the labeled and
nonplane counterpart of classical Motzkin trees. For a comprehensive introduction into recursive
structures like Motzkin trees and also labeled and nonplane combinatorial structures see
\cite{AnaCombi}.}, \emph{i.e.}, a tree consisting of leaves (zero children), unary vertices (one
child) and binary vertices (two children). All edges in this Motzkin tree are directed away from
the root. In particular, it is a Motzkin tree with exactly $2k$ unary vertices, where $k$ of them
are the starting points of the removed edges, the other $k$ their end points (note that here the
tree-child property was used).

\begin{figure}[ht]
\begin{center}
\begin{tikzpicture}[scale=0.7]
\hspace*{-3cm}
\draw (9cm,0cm) node[inner sep=1.5pt,circle,draw,fill] (1) {};
\draw (8cm,-1cm) node[inner sep=1.5pt,circle,draw,fill] (2) {};
\draw (10cm,-1cm) node[color=green!70,inner sep=1.5pt,circle,draw,fill] (3) {};
\draw (7cm,-2cm) node[inner sep=1.5pt,circle,draw,fill] (4) {};
\draw (9cm,-2cm) node[color=red!70,inner sep=1.5pt,circle,draw,fill] (5) {};
\draw (11cm,-2cm) node[inner sep=1.5pt,circle,draw,fill] (6) {};
\draw (6cm,-3cm) node[inner sep=1.5pt,circle,draw,fill] (7) {};
\draw (7cm,-3.5cm) node[inner sep=1.5pt,circle,draw,fill] (9) {};
\draw (9cm,-3.5cm) node[color=green!70,inner sep=1.5pt,circle,draw,fill] (10) {};
\draw (6cm,-4.5cm) node[inner sep=1.5pt,circle,draw,fill] (12) {};
\draw (8cm,-4.5cm) node[color=red!70,inner sep=1.5pt,circle,draw,fill] (13) {};
\draw (10cm,-4.5cm) node[inner sep=1.5pt,circle,draw,fill] (14) {};
\draw (8cm,-5.7cm) node[inner sep=1.5pt,circle,draw,fill] (16) {};
\draw (9.27cm,-1.33cm) node[inner sep=1.5pt,circle] (21) {};
\draw (9.67cm,-1.73cm) node[inner sep=1.5pt,circle] (22) {};
\draw (9.33cm,-1.27cm) node[inner sep=1.5pt,circle] (23) {};
\draw (9.73cm,-1.67cm) node[inner sep=1.5pt,circle] (24) {};
\draw (8.27cm,-3.83cm) node[inner sep=1.5pt,circle] (25) {};
\draw (8.67cm,-4.23cm) node[inner sep=1.5pt,circle] (26) {};
\draw (8.33cm,-3.77cm) node[inner sep=1.5pt,circle] (27) {};
\draw (8.73cm,-4.17cm) node[inner sep=1.5pt,circle] (28) {};

\draw (1)--(2);\draw (1)--(3);\draw (2)--(4);\draw (2)--(5);
\draw (3)--(5);\draw (3)--(6);\draw (4)--(7);
\draw (4)--(9);\draw (5)--(10);\draw (9)--(12);
\draw (9)--(13);\draw (10)--(13);\draw (10)--(14);
\draw (13)--(16);\draw (21)--(22);\draw (23)--(24);
\draw (25)--(26);\draw (27)--(28);

\draw (13cm,-2.85cm) node[inner sep=1.5pt,circle] (29) {$\longrightarrow$};

\hspace*{6.4cm}
\draw (9cm,0cm) node[inner sep=1.5pt,circle,draw,fill] (1) {};
\draw (8cm,-1cm) node[inner sep=1.5pt,circle,draw,fill] (2) {};
\draw (10cm,-1cm) node[color=green!70,inner sep=1.5pt,circle,draw,fill] (3) {};
\draw (7cm,-2cm) node[inner sep=1.5pt,circle,draw,fill] (4) {};
\draw (9cm,-2cm) node[color=red!70,inner sep=1.5pt,circle,draw,fill] (5) {};
\draw (11cm,-2cm) node[inner sep=1.5pt,circle,draw,fill] (6) {};
\draw (6cm,-3cm) node[inner sep=1.5pt,circle,draw,fill] (7) {};
\draw (7cm,-3.5cm) node[inner sep=1.5pt,circle,draw,fill] (9) {};
\draw (9cm,-3.5cm) node[color=green!70,inner sep=1.5pt,circle,draw,fill] (10) {};
\draw (6cm,-4.5cm) node[inner sep=1.5pt,circle,draw,fill] (12) {};
\draw (8cm,-4.5cm) node[color=red!70,inner sep=1.5pt,circle,draw,fill] (13) {};
\draw (10cm,-4.5cm) node[inner sep=1.5pt,circle,draw,fill] (14) {};
\draw (8cm,-5.7cm) node[inner sep=1.5pt,circle,draw,fill] (16) {};

\draw (1)--(2);\draw (1)--(3);\draw (2)--(4);\draw (2)--(5);
\draw (3)--(6);\draw (4)--(7);
\draw (4)--(9);\draw (5)--(10);\draw (9)--(12);
\draw (9)--(13);\draw (10)--(14);\draw (13)--(16);

\draw (13cm,-2.85cm) node[inner sep=1.5pt,circle] (29) {$\longrightarrow$};

\draw (15.8cm,-2.35cm) node[inner sep=1.5pt,circle,draw,fill] (30) {};
\draw (14.8cm,-3.35cm) node[color=green!70,inner sep=1.5pt,circle,draw,fill] (31) {};
\draw (16.8cm,-3.35cm) node[color=green!70,inner sep=1.5pt,circle,draw,fill] (32) {};

\draw (30)--(31);\draw (30)--(32);
\end{tikzpicture}
\end{center}
\caption{A normal network with colored Motzkin skeleton and coresponding sparsened skeleton. Note that there are three more possible colored Motzkin skeletons which one can obtain from the same network and that all but one yield the same sparsened skeleton.}
\label{decomp_example}
\end{figure}
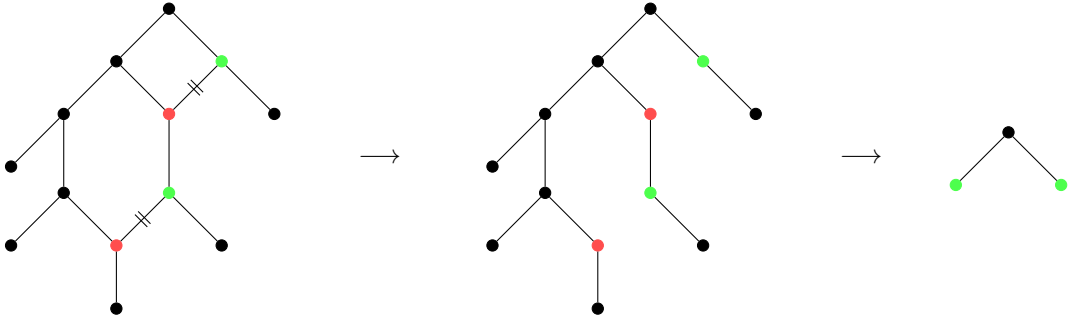

Now consider the following procedure (see Figure~\ref{decomp_example} for an illustration):
Start with a Motzkin tree $M$ with exactly $2k$ unary vertices
and $n$ vertices in total. Then add edges such that {\it (i)} each edge connects two unary vertices,
{\it (ii)} no  two of the added edges have a vertex in common,
and {\it (iii)} the resulting graph is a normal network $N$. Finally, color the start
vertices of the added edges green and their end vertices red.
We say then that $M$ (keeping the colors from the above generation of $N$, but not the edges)
is a {\em colored Motzkin skeleton} (or simply {\em Motzkin skeleton}) of $N$.
In this way all normal networks with $n$ vertices
are generated and each of them exactly $2^k$ times, since every network $N$ with $ k $ reticulation vertex has exactly $2^k$
different Motzkin skeletons.

In order to set up generating functions for phylogenetic networks, we will construct them as
follows: for a given network $N$ first pick one of its $2^k$ possible Motzkin skeletons. Then, look for the minimal
subtree $T$ which contains all green vertices. This tree contains all the green vertices as well
as all last common ancestors\footnote{Note that we use the name which is common in the
combinatorial literature. In the phylogenetics literature this is usually called \emph{most recent
common ancestor}.} of any two green
vertices. These particular vertices form a tree whose edges are paths in $T$. Contract each of these paths to one single edge.
The resulting tree, which is again a Motzkin tree, is called the {\em sparsened skeleton} of $N$. The structure of this tree tells us how the green vertices are distributed within $N$ (again See Figure~\ref{decomp_example}).

In order to construct networks with $k$ reticulation vertices, we start with a sparsened skeleton
having $k$ green vertices. Then we replace all edges by paths that are made of red vertices or binary vertices with a
Motzkin tree (whose unary vertices are all colored red) as second child and add a path of the same type on top of the root of the sparsened
skeleton. Moreover, we attach a Motzkin tree (again with all unary vertices colored red) to each leaf of the
sparsened skeleton such that this new subtree is linked to the
(now former) leaf by an edge (for normal networks, this tree can be a binary tree).
Do all of the above in such a way that the new structure has $k$ red vertices altogether.
What we obtain so far is a Motzkin skeleton of a phylogenetic network. Finally, add
edges connecting the green vertices to the red ones in such a way that the corresponding mapping
is bijective and that the normality condition for phylogenetic networks is respected.

Let us set up the exponential generating function for Motzkin trees which appear in the above
construction. This means that the tree-child condition for networks has to be respected, but
the number of unary vertices need not  be even. After all, the unary vertices in those trees will be
the red vertices of our network.

Denote by $M_{\ell,n}$ the number of all vertex-labeled Motzkin trees with $n$ vertices
and $\ell$ unary vertices (all colored red) that respect the tree-child condition for networks, which means that the child of a unary vertex
cannot be a unary vertex and each binary vertex has at least one child which is not a unary vertex.
Let $\mathcal M$ denote the set of all these Motzkin trees.
The exponential generating function associated to $\mathcal M$ is
\[
M(z,y)=\sum_{n\ge 1}\sum_{\ell\ge 0} M_{\ell,n} y^\ell \frac{z^n}{n!}.
\]
Furthermore, let $M_u(z,y)$ and $M_b(z,y)$ denote the generating function associated to all
Motzkin trees in $\M$ whose root is a unary vertex and a binary vertex, respectively. Then
\begin{equation*}
M_u(z,y)=zy(z+M_b(z,y))
\end{equation*}
since a unary vertex cannot have a unary child. In a Motzkin tree with a binary root, the root may have two children being either a
leaf or a binary vertex, or one of the children is a unary vertex and the other either a leaf or a
binary vertex. This yields
\begin{equation*}
M_b(z,y)=\frac z2( (z+M_b(z,y))^2+2zy(z+M_b(z,y))^2).
\end{equation*}
Solving gives
\[
M_b(z,y)=\frac{1-\sqrt{1-2z^2-4yz^3}}{z(1+2yz)}-z
\]
and
\begin{equation}\label{Mu}
M_u(z,y)=y\frac{1-\sqrt{1-2z^2-4yz^3}}{1+2yz}
\end{equation}
and thus
\begin{equation} \label{fctM}
M(z,y)=z+M_u(z,y)+M_b(z,y)=\frac{\(1+yz\)\(1-\sqrt{1-2z^2-4yz^3}\,\)}{z(1+2yz)}.
\end{equation}
The first few coefficients can be seen from
\[
z+yz^2+\frac12 z^3+\frac32 yz^4+\(y^2+\frac12\)z^5+\frac52
yz^6+\(4y^2+\frac58\)z^7+\(2y^3+\frac{35}8y\)z^8+\cdots.
\]

\section{Counting Vertex-Labeled Normal Networks}\label{vl-normal-net}

In this section, we will count (vertex-labeled) normal networks with a fixed number $k$ of reticulation vertices. We will start with the cases $k=1,2,3$ which will be discussed in the next three subsections and for which we will derive asymptotic expansions up to the second order term (in fact, our method allows one to obtain full asymptotic expansions as well). From these three cases, we will observe a general pattern which will be proved in the fourth subsection.

\subsection{Normal networks with one reticulation vertex}

In this subsection we will determine the asymptotic number of normal networks with one reticulation
vertex and then discuss their relationship to unicyclic networks that were studied in \cite{SeSt06}.

\subsubsection{Counting}

In order to count normal networks with only one reticulation vertex we use Motzkin trees from the class ${\mathcal M}$, which have generating function
\eqref{fctM}, and (sparsened) skeletons, as described in the previous section: We delete
one of the two incoming edges of the reticulation vertex which then gives a unary-binary tree
satisfying the tree-child property with exactly two unary vertices. Conversely, we can start with the
general tree or even the sparsened skeleton (which only consists of one vertex) and then construct the network from this.

\bp\label{exp-k=1}
The exponential generating function for vertex-labeled normal networks with
one reticulation vertex is
\begin{equation} \label{fctN1}
N_1(z)
=\frac{z\(1-\sqrt{1-2z^2}\,\)^3}{2(1-2z^2)^{3/2}}=z\frac{\tilde{a}_1(z^2)-\tilde{b}_1(z^2)\sqrt{1-2z^2}}{(1-2z^2)^{3/2}},
\end{equation}
where
\[
\tilde{a}_1(z)=2-3z\qquad\text{and}\qquad\tilde{b}_1(z)=2-z.
\]
\ep

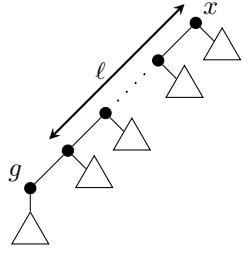
\begin{figure}
	\begin{center}
		\begin{tikzpicture}[
		scale=0.5,
		level/.style={thick},
		virtual/.style={thick,densely dashed},
		trans/.style={thick,<->,shorten >=2pt,shorten <=2pt,>=stealth},
		classical/.style={thin,double,<->,shorten >=4pt,shorten <=4pt,>=stealth}
		]		
		
		\draw (0cm,0cm) node[inner sep=1.5pt,circle,draw,fill] (1) {};
		\draw (1cm,1cm) node[inner sep=1.5pt,circle,draw,fill] (2) {};
		\draw (2cm,2cm) node[inner sep=1.5pt,circle,draw,fill] (3) {};
		\draw (3.4cm,3.4cm) node[inner sep=1.5pt,circle,draw,fill] (4) {};
		\draw (4.4cm,4.4cm) node[inner sep=1.5pt,circle,draw,fill] (5) {};
		\draw (2.2cm,2.2cm) node[inner sep=1.5pt,circle] (11) {};
		\draw (3.2cm,3.2cm) node[inner sep=1.5pt,circle] (12) {};
		\draw (4.8cm,4.8cm) node[inner sep=1.5pt,circle] (14) {$x$};
		\draw (-0.4cm,0.4cm) node[inner sep=1.5pt,circle] (13) {$g$};
		\draw (0cm,-1.2cm) node[regular polygon,regular polygon sides=3,draw,inner sep=0.1cm] (6) {};
		\draw (1.7cm,0.3cm) node[regular polygon,regular polygon sides=3,draw,inner sep=0.1cm] (7) {};
		\draw (2.7cm,1.3cm) node[regular polygon,regular polygon sides=3,draw,inner sep=0.1cm] (8) {};
		\draw (4.1cm,2.7cm) node[regular polygon,regular polygon sides=3,draw,inner sep=0.1cm] (9) {};
		\draw (5.1cm,3.7cm) node[regular polygon,regular polygon sides=3,draw,inner sep=0.1cm] (10) {};
		
		\draw (1)--(2);\draw (2)--(3);\draw (4)--(5);\draw (1)--(6);
		\draw (2)--(7);\draw (3)--(8);\draw (4)--(9);\draw (5)--(10);
		\draw[loosely dotted,line width=0.7pt] (11)--(12);
		\draw[trans] (0.4cm,1.2cm) -- (4.2cm,5cm) node[midway,left] {\El{}};	
		\end{tikzpicture}
	\end{center}
	\caption{The structure of Motzkin skeletons of networks with one
reticulation vertex. It originates from a sparsened skeleton which consists of only one green vertex. It has one
green vertex, denoted by $g$, and one red vertex which is hidden within the forest made of the
triangles in the picture, which are attached to $g$ and all the vertices on the path of length $\ell$.
Note that the position of the red vertex in this forest is restricted by the normality condition.} \label{fig_skeleton}
\end{figure}

\bpf
As already mentioned, we start with the general tree as depicted in Figure~\ref{fig_skeleton},
which arises from the sparsened skeleton, \emph{i.e.}, the tree consisting of a single green vertex $g$ as follows: we add a sequence of trees on top of $g$ which consist of a root to which a tree in ${\mathcal M}$ is attached. Moreover, we attach also a tree from ${\mathcal M}$ to $g$ as a subtree.

Next, in order to obtain all normal networks arising from these Motzkin skeletons, we have to add an edge starting from $g$ and
ending at the red vertex. Note that for a normal network, this edge is neither allowed to
point to a vertex on the path from $g$ to the root (since the network must be a DAG), nor to the
root of one of the trees which are connected to the vertices on the path from $g$ to the root (since
this violates the normality condition) nor to any vertex in the subtree of $g$ (since this again
violates the normality condition). Overall, the red vertex must be contained in the forest attached to the path from
$g$ to the root, but not in the tree attached to $g$.
Moreover, note that since there is only one red vertex, the requirement that trees in this forest satisfy the tree-child property could actually be dropped.

The networks arising from these skeletons can therefore be specified as a tree without red vertices (the one attached to $g$) and a sequence of structures of the form ``vertex plus Motzkin tree
with non-unary root'' (\emph{cf}. Figure~\ref{fig_skeleton}). In terms of generating functions this gives

\[
N_1(z)=
\frac12 \frac{\partial}{\partial y} \frac{zM(z,0)}{1-z\tilde M(z,y)}\Big\vert_{y=0}
\]
where
\begin{align}\label{tree-without-root}
\tilde{M}(z,y)&=z+M_b(z,y)=M(z,y)-zy(z+M_b(z,y)).
\end{align}
The factor $1/2$ makes up for the fact that each network
is counted exactly twice by the above procedure. Evaluating this and writing $M_y$ for the partial derivative of $M$ (w.r.t. $y$) yields
\begin{align}
N_1(z)&=\frac{z}{2}M(z,0)(M_y(z,0)-z^2-zM_b(z,0))\sum_{\ell\geq 1}\ell z^{\ell}M(z,0)^{\ell-1} \nonumber\\
&=\frac{z^2M(z,0)(M_y(z,0)-z^2-zM_b(z,0))}{2(1-zM(z,0))^2}. \nonumber
\end{align}
Now, by using
\begin{equation}\label{exp-M-y-b}
M(z,0)=\frac{1-\sqrt{1-2z^2}}{z},\qquad M_{y}(z,0)=\frac{1}{\sqrt{1-2z^2}}-1,\qquad M_b(z,0)=\frac{1-\sqrt{1-2z^2}}{z}-z,
\end{equation}
we obtain \eqref{fctN1}.
\epf

From this result we can now easily obtain the asymptotic number of normal networks (see the appendix for numerical data).

\begin{cor}\label{asymp-k=1}
Let $N_{1,n}$ denote the number of vertex-labeled normal networks with $n$ vertices
and one reticulation vertex. If $n$ is even then $N_{1,n}$ is zero, otherwise
\[
N_{1,n}=n![z^n]N_1(z)=\(\frac{\sqrt2}{e}\)^n n^{n+1}\(\frac{\sqrt{2}}{2}-\frac{3\sqrt{\pi}}{2}\cdot\frac{1}{\sqrt{n}}+{\mathcal O}\left(\frac{1}{n}\right)\),
\]
as $\nti$.
\end{cor}
\bpf
The function \eqref{fctN1} has two dominant singularities, namely at $\pm1/\sqrt{2}$, with
singular expansions
\[
N_1(z)\stackrel{z\rightarrow\pm1/\sqrt{2}}{\sim}\pm\frac{1}{8(1\mp\sqrt{2}z)^{3/2}}\mp\frac{3\sqrt{2}}{8(1\mp\sqrt{2}z)}+{\mathcal O}\left(\frac{1}{\sqrt{1\mp\sqrt{2}z}}\right).
\]
Applying a transfer lemma (see \cite{FO90, AnaCombi})
for these two singularities and using Stirling's formula completes the proof.
\epf

\brem
Note that the periodicity is not surprising since, as mentioned in the introduction, phylogenetic networks always have an odd number of
vertices.
\erem

\subsubsection{Relationship to unicyclic networks}\label{rel-uni}

In \cite{SeSt06}, the authors counted unicyclic networks which are (vertex-labeled or leaf-labeled)
pointed\footnote{Pointed means that an edge is chosen to which a vertex $v$ is
attached (with an edge, of course). The chosen edge itself is thus split into two edges and the
point where the new edge is attached becomes a further new vertex. The vertex $v$ is then the
root vertex of the network} graphs with only one
cycle to which complete binary trees are attached. The enumeration was done only for
leaf-labeled networks there.

On the other hand, phylogenetic networks with exactly one reticulation vertex are the
same as unicyclic networks, if one disregards the direction of the
edges.\footnote{Combinatorially, there is no big difference between rooted and pointed, since we
can always drop the attached vertex and edge in the latter case and direct all edges. Thus, if one can solve
the counting problem for a subclass of rooted networks also the corresponding counting problem for pointed graphs can be solved.}

So, another way of counting normal networks with exactly one reticulation vertex is by using a
modification of the approach of \cite{SeSt06}: either the root is in a cycle in which case one of
the vertices on this cycle except the root and its two neighbours are the reticulation vertex and to
each vertex may be attached a complete binary rooted tree or the root is not in the cycle in which
case exactly one subtree contains the cycle. This translates into
\[
N_1(z)=zM(z,0)N_1(z)+\frac{1}{2}\sum_{\ell\geq 3}(\ell-2)z^{\ell+1}M(z,0)^{\ell}.
\]
Solving this equation gives
\begin{equation*}
N_1(z)=\frac{\sum_{\ell\geq
3}(\ell-2)z^{\ell+1}M(z,0)^{\ell}}{2(1-zM(z,0))}=\frac{z^4M(z,0)^3}{2(1-zM(z,0))^3}.
\end{equation*}
Plugging (\ref{exp-M-y-b}) into this reveals
\[
N_1(z)=\frac{z(1-\sqrt{1-2z^2})^3}{2(1-2z^2)^{3/2}}
\]
as it must be.

\subsection{Normal networks with two reticulation vertices}\label{nnt}

For this case, we use two variables $y_1,y_2$ to express the possible pointings of the two green vertices of the Motzkin skeletons.
Furthermore, we have now more complicated paths (and attached trees) which replace the edges of the sparsened skeleton and thus we first set up the generating function corresponding to theses paths. To govern the situation where an edge from one of the two green vertices must not point into a
certain subtree or to a particular vertex, we distinguish several types of unary vertices, which are the red vertices of our construction.

To simplify the explanation, let us use the following conventions: If the root of a Motzkin tree is a unary vertex (so, a red vertex) we call the tree a
\emph{red tree}, otherwise a \emph{white tree}. Note that the class of white trees has generating function $\tilde M(z,y)$ given in \eqref{tree-without-root}, whereas the class of red trees has generating function $M_u(z,y)$, see \eqref{Mu}.

The structure we will need is a class $\mathcal P$ of paths which serve as the essential building blocks for Motzkin skeletons. In this class the rules for pointing to particular red vertices differ, depending on whether (i) the red vertex lies on the path itself, but is not the very first vertex there, or is the root of one of the (red) subtrees attached to the vertices of the path, (ii) it is one of the non-root vertices of one of the attached subtrees or (iii) the red vertex is the first vertex of the path. To distinguish these three classes of red vertices, we will mark the red vertices of type (i) with the variable $y$, those of type
(ii) with $\tilde y$ and the vertex of type (iii) with $\hat y$.

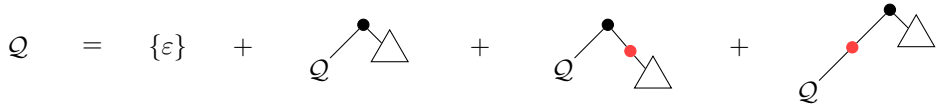
\begin{figure}[h]
\begin{center}
 \begin{tikzpicture}[
                scale=0.5,
                level/.style={thick},
                virtual/.style={thick,densely dashed},
                trans/.style={thick,<->,shorten >=2pt,shorten <=2pt,>=stealth},
                classical/.style={thin,double,<->,shorten >=4pt,shorten <=4pt,>=stealth}
                ]

                \draw (0cm,0.8cm) node[inner sep=1.5pt,circle] (2) {};
        \draw (-0.2cm,0.6cm) node[inner sep=1.5pt,circle] (1) {$\mathcal Q$};
                \draw (1cm,1.8cm) node[inner sep=1.5pt,circle,draw,fill] (3) {};
                \draw (1.7cm,1.1cm) node[regular polygon,regular polygon sides=3,draw,inner sep=0.1cm] (4) {};

                \draw (2)--(3);\draw (3)--(4);

        \draw (6.5cm,0.8cm) node[inner sep=1.5pt,circle] (2) {};
        \draw (6.3cm,0.6cm) node[inner sep=1.5pt,circle] (1) {$\mathcal Q$};
                \draw (7.5cm,1.8cm) node[inner sep=1.5pt,circle,draw,fill] (3) {};
                \draw (8.7cm,0.4cm) node[regular polygon,regular polygon sides=3,draw,inner sep=0.1cm] (5) {};
                \draw (8.1cm,1.1cm) node[inner sep=1.5pt,circle,draw,fill,red!75] (4) {};

                \draw (2)--(3);\draw (3)--(4);\draw (4)--(5);

            \draw (13cm,0.2cm) node[inner sep=1.5pt,circle] (2) {};
        \draw (12.8cm,0cm) node[inner sep=1.5pt,circle] (1) {$\mathcal Q$};
                \draw (15cm,2.2cm) node[inner sep=1.5pt,circle,draw,fill] (3) {};
                \draw (15.7cm,1.5cm) node[regular polygon,regular polygon sides=3,draw,inner sep=0.1cm] (4) {};
        \draw (14cm,1.2cm) node[inner sep=1.5pt,circle,draw,fill,red!75] (5) {};

        \draw (2)--(5);\draw (5)--(3);\draw (3)--(4);

        \draw (4.1cm,1.2cm) node[inner sep=1.5pt,circle] (1) {$+$};
        \draw (11cm,1.2cm) node[inner sep=1.5pt,circle] (1) {$+$};
        \draw (-2.2cm,1.2cm) node[inner sep=1.5pt,circle] (1) {$+$};
        \draw (-4.2cm,1.2cm) node[inner sep=1.5pt,circle] (1) {$\{\varepsilon\}$};
        \draw (-6.2cm,1.2cm) node[inner sep=1.5pt,circle] (1) {$=$};
        \draw (-8.2cm,1.2cm) node[inner sep=1.5pt,circle] (1) {$\mathcal Q$};
                \end{tikzpicture}
\end{center}	
\caption{The specification of the class $\mathcal Q$. In this picture, the paths are drawn such that they are going from upper right to lower left. The triangles represent the trees attached to the path. These are white trees, \emph{i.e.}, trees which do not have a
unary root. The variable $y$ marks the red vertices that are shown in the figure. Others may be hidden in the white trees and are marked by $\tilde y$.
The last part of the specification guarantees that there are no consecutive red vertices on the path.}	
\label{QGF}
\end{figure}

Moreover, we have to respect the tree-child condition. Normality does not play a role on this level, it actually only causes the need for the third class of red vertices. The tree-child condition implies that the successor of a red vertex on the path itself must not be red. Moreover, if the tree attached to some vertex $x$ is a red tree, then the successor of $x$ on the path must not be a red vertex. This
gives rise to the a combinatorial specification. Take a set of three possible atomic items: a vertex with a white tree, a vertex with a red tree (which is itself a red vertex having a white tree as subtree), and a vertex having a red vertex and a white tree as subtrees. Then a path in $\mathcal P$ is either (a) a sequence made of these atomic items or (b) a red vertex followed by a sequence of type (a). More formally, let $\tilde{\mathcal M}$ denote the class of white trees, $\circ$ denote a binary vertex and $\bullet$ denote a red (unary) vertex. We write $\{x\}\times S \times T$ if $x$ is a vertex having subtrees $S$ and $T$, where $T$ is omitted if $x$ is a red vertex and the edge $x$ --- $S$ is an edge of the path. Then we consider a class $\mathcal Q$ which contains all path in ${\mathcal P}$ of type (a) above. The specification of this class is
\begin{equation}\label{specQ}
\mathcal Q=\{\eps\}\cup \{\circ\}\times\mathcal Q\times \tilde{\mathcal M} \cup \{\circ\}\times\mathcal Q\times (\{\bullet\}\times\tilde{\mathcal M})
\cup \{\circ\}\times(\{\bullet\}\times \mathcal Q) \times \tilde{\mathcal M},
\end{equation}
where $\varepsilon$ denotes the empty tree; see Figure~\ref{QGF}. Since a path in ${\mathcal P}$ may also start with a red vertex, which then belongs to the third class of red vertices, we specify $\mathcal P$ as
\begin{equation}\label{specP}
\mathcal P=\mathcal Q \cup \{\bullet\}\times \mathcal Q.
\end{equation}
This leads to the generating function
\begin{equation}\label{PGF}
P(z,y,\tilde y,\hat y)=\frac{1+z\hat y}{1-(z+2z^2y)\tilde M(z,\tilde y)}
\end{equation}
after all.

Let us summarize what we just defined. In our analysis the variables $y$, $\tilde y$, $\hat y$ will be replaced by a
sum of variables $y_i$ where the presence of a particular $y_i$ indicates that the corresponding $g_i$ is allowed to point,
its absence that pointing is forbidden. In particular, $y$
represent the permission to point to vertices of the path (except its first vertex) as well as to the roots of the trees attached to the path. The variable
$\tilde y$ describes the permission to point to non-root vertices of these trees
and $\hat y$ allows pointing to the first vertex of the path.

Now we are ready to state the following result.

\bp\label{exp-k=2}
The exponential generating function for vertex-labeled normal networks with
two reticulation vertices is
\begin{equation} \label{fctN2}
N_2(z)=z\frac{\tilde{a}_2(z^2)-\tilde{b}_2(z^2)\sqrt{1-2z^2}}{(1-2z^2)^{7/2}},
\end{equation}
where
\begin{equation} \label{polyAB}
\tilde{a}_2(z)=6z^4-\frac{87}{2}z^3+30z^2-3z
\quad\text{and}\quad
\tilde{b}_2(z)=-18z^3+27z^2-3z
\end{equation}
\ep

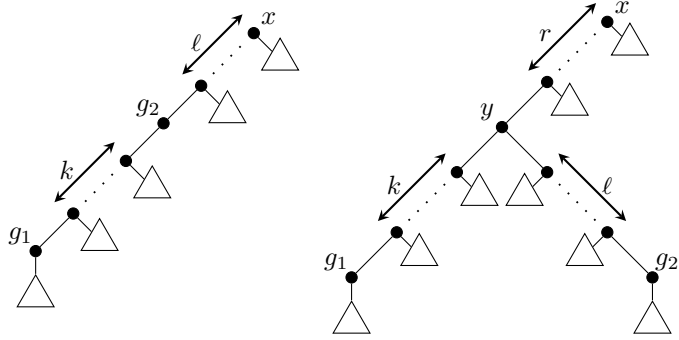
\begin{figure}
	\begin{center}
		\begin{tikzpicture}[
		scale=0.5,
		level/.style={thick},
		virtual/.style={thick,densely dashed},
		trans/.style={thick,<->,shorten >=2pt,shorten <=2pt,>=stealth},
		classical/.style={thin,double,<->,shorten >=4pt,shorten <=4pt,>=stealth}
		]
		
        \draw (-11.2cm,-2.1cm) node[inner sep=1.5pt,circle,draw,fill] (1) {};
        \draw (-11.6cm,-1.7cm) node[inner sep=1.5pt,circle] (13) {$g_{1}$};
		\draw (-10.2cm,-1.1cm) node[inner sep=1.5pt,circle,draw,fill] (2) {};
		\draw (-8.8cm,0.3cm) node[inner sep=1.5pt,circle,draw,fill] (3) {};
		\draw (-7.8cm,1.3cm) node[inner sep=1.5pt,circle,draw,fill] (4) {};
		\draw (-8.2cm,1.7cm) node[inner sep=1.5pt,circle] (14) {$g_{2}$};
		\draw (-6.8cm,2.3cm) node[inner sep=1.5pt,circle,draw,fill] (5) {};
		\draw (-5.4cm,3.7cm) node[inner sep=1.5pt,circle,draw,fill] (6) {};
		\draw (-5cm,4.1cm) node[inner sep=1.5pt,circle] (15) {$x$};
        \draw (-11.2cm,-3.3cm) node[regular polygon,regular polygon sides=3,draw,inner sep=0.1cm] (7) {};
		\draw (-9.5cm,-1.8cm) node[regular polygon,regular polygon sides=3,draw,inner sep=0.1cm] (8) {};
		\draw (-8.1cm,-0.4cm) node[regular polygon,regular polygon sides=3,draw,inner sep=0.1cm] (9) {};
		\draw (-6.1cm,1.6cm) node[regular polygon,regular polygon sides=3,draw,inner sep=0.1cm] (11) {};
        \draw (-4.7cm,3cm) node[regular polygon,regular polygon sides=3,draw,inner sep=0.1cm] (12) {};
        \draw (-10cm,-0.9cm) node[inner sep=1.5pt,circle] (13) {};
        \draw (-9cm,0.1cm) node[inner sep=1.5pt,circle] (14) {};
        \draw (-6.6cm,2.5cm) node[inner sep=1.5pt,circle] (15) {};
        \draw (-5.6cm,3.5cm) node[inner sep=1.5pt,circle] (16) {};
		
		\draw (1)--(2);\draw (3)--(4);\draw (4)--(5);
        \draw (1)--(7);\draw (2)--(8);\draw (3)--(9);
        \draw (5)--(11);\draw (6)--(12);
		\draw[loosely dotted,line width=0.7pt] (13)--(14);
		\draw[loosely dotted,line width=0.7pt] (15)--(16);
		\draw[trans] (-10.8cm,-0.9cm) -- (-9cm,0.9cm) node[midway,left] {\Ke{}};
		\draw[trans] (-7.4cm,2.5cm) -- (-5.6cm,4.3cm) node[midway,left] {\El{}};

        scale=0.6
        \draw (0cm,0cm) node[inner sep=1.5pt,circle,draw,fill] (1) {};
        \draw (1.2cm,1.2cm) node[inner sep=1.5pt,circle,draw,fill] (2) {};
        \draw (0.8cm,1.6cm) node[inner sep=1.5pt,circle] (26) {$y$};
        \draw (2.4cm,0cm) node[inner sep=1.5pt,circle,draw,fill] (3) {};
        \draw (0.6cm,-0.6cm) node[regular polygon,regular polygon sides=3,draw,inner sep=0.1cm] (4) {};
        \draw (1.8cm,-0.6cm) node[regular polygon,regular polygon sides=3,draw,inner sep=0.1cm] (5) {};
        \draw (-1.6cm,-1.6cm) node[inner sep=1.5pt,circle,draw,fill] (6) {};
        \draw (4cm,-1.6cm) node[inner sep=1.5pt,circle,draw,fill] (7) {};
        \draw (-2.8cm,-2.8cm) node[inner sep=1.5pt,circle,draw,fill] (8) {};
        \draw (-3.2cm,-2.4cm) node[inner sep=1.5pt,circle] (24) {$g_1$};
        \draw (5.2cm,-2.8cm) node[inner sep=1.5pt,circle,draw,fill] (9) {};
        \draw (5.6cm,-2.4cm) node[inner sep=1.5pt,circle] (25) {$g_2$};
        \draw (-1cm,-2.2cm) node[regular polygon,regular polygon sides=3,draw,inner sep=0.1cm] (10) {};
        \draw (3.4cm,-2.2cm) node[regular polygon,regular polygon sides=3,draw,inner sep=0.1cm] (11) {};
        \draw (-2.8cm,-4cm) node[regular polygon,regular polygon sides=3,draw,inner sep=0.1cm] (12) {};
        \draw (5.2cm,-4cm) node[regular polygon,regular polygon sides=3,draw,inner sep=0.1cm] (13) {};
        \draw (2.4cm,2.4cm) node[inner sep=1.5pt,circle,draw,fill] (14) {};
        \draw (4cm,4cm) node[inner sep=1.5pt,circle,draw,fill] (15) {};
        \draw (4.4cm,4.4cm) node[inner sep=1.5pt,circle] (27) {$x$};
        \draw (3cm,1.8cm) node[regular polygon,regular polygon sides=3,draw,inner sep=0.1cm] (16) {};
        \draw (4.6cm,3.4cm) node[regular polygon,regular polygon sides=3,draw,inner sep=0.1cm] (17) {};
        \draw (-1.4cm,-1.4cm) node[inner sep=1.5pt,circle] (18) {};
        \draw (-0.2cm,-0.2cm) node[inner sep=1.5pt,circle] (19) {};
        \draw (3.8cm,-1.4cm) node[inner sep=1.5pt,circle] (20) {};
        \draw (2.6cm,-0.2cm) node[inner sep=1.5pt,circle] (21) {};
        \draw (2.6cm,2.6cm) node[inner sep=1.5pt,circle] (22) {};
        \draw (3.8cm,3.8cm) node[inner sep=1.5pt,circle] (23) {};

        \draw (8)--(12);\draw (8)--(6);\draw (6)--(10);\draw (1)--(4);\draw (1)--(2);
        \draw (2)--(3);\draw (3)--(5);\draw (7)--(11);\draw (7)--(9);\draw (9)--(13);
        \draw (2)--(14);\draw (14)--(16);\draw (15)--(17);
        \draw [loosely dotted,line width=0.7pt] (18)--(19);
        \draw [loosely dotted,line width=0.7pt] (20)--(21);
        \draw [loosely dotted,line width=0.7pt] (22)--(23);
        \draw[trans] (-2.2cm,-1.4cm) -- (-0.2cm,0.6cm) node[midway,left] {\Ke{}};
		\draw[trans] (2.6cm,0.6cm) -- (4.6cm,-1.4cm) node[midway,right] {\El{}};
        \draw[trans] (1.8cm,2.6cm) -- (3.8cm,4.6cm) node[midway,left] {\Rl{}};
		\end{tikzpicture}
	\end{center}
	\caption{The possible structures of Motzkin skeletons of networks with
two reticulation vertices. These originate from the two possible sparsened skeletons made
of two green vertices: The path of length one, which gives rise to the left Motzkin skeleton, and
the cherry leading to the right Motzkin skeleton.
\newline {\bf Note:} In this figure (as well as in all the subsequent figures of this paper)
the triangles are placeholders for trees which may but need not necessarily be white trees (see beginning of Section~\ref{nnt}).
The class they belong to depends on their position within the normal network.} \label{fig_skeleton_2}
\end{figure}

\bpf
Note that, in the current situation, there are only two possible sparsened skeletons: either a path of length one (with both vertices green) or a cherry (with both leaves being green vertices). From this, one builds two possible types of Motzkin skeletons that
are depicted in Figure~\ref{fig_skeleton_2}.

For the first type (see Figure~\ref{fig_skeleton_2}, left), we
have to complete the Motzkin skeletons by adding
two egdes having start vertex $g_1$ and $g_2$, respectively. The one starting from $g_1$ may point
to any non-root vertex within the subtrees that are attached to the skeleton's spine (\emph{i.e.}, the paths $k$ and $\ell$ and $g_2$). By normality, it can neither point to the root of one of those subtrees nor to a vertex in the subtree attached to $g_1$ itself, but $g_1$ does not
belong to what we called the spine anyway. Similarly, the edge starting at $g_2$ may point to any non-root vertex in the subtrees attached to the path $\ell$, the path from
the root to the parent of $g_2$. Thus the generating function of the subtrees attached to the vertices of $\ell$ is $\tilde M(z,y_1+y_2)$, that of the subtrees  attached to the vertices of $k$ is $\tilde M(z,y_1)$. The tree attached to $g_1$ corresponds to $M(z,0)$ since it must not contain any red vertices.
Finally, note that we have to point at two red vertices,
one targeted by $g_1$ and one targeted by $g_2$. Pointing (and not counting it any more as red vertex) corresponds to differentiation in the world of generating functions. Since we do not want any other red vertices to be present, we set $y_1=y_2=0$ after the differentiations. After all, we get the generating function
\begin{align*}
N_{2,1}(z)
&=\partial_{y_1} \partial_{y_2}\frac{z^2M(z,0)}{(1-z\tilde M(z,y_1))(1-z\tilde
M(z,y_1+y_2))}\Big\vert_{y_1=0,y_2=0}.
\end{align*}

For the second type (see Figure~\ref{fig_skeleton_2}, right), none of the two green vertices $g_1$ and $g_2$ is the ancestor
of the other, they have a common ancestor $y$. Moreover, there is a path on top of $y$ connecting $y$ with the root
of the network, called $r$. Also, there are paths from $y$ to $g_1$ and $g_2$, namely $k$
and $\ell$, respectively. To each of the vertices of $k$, $\ell$ and $r$ as well as to the two green vertices
a tree from ${\mathcal M}$ is attached.

In order to meet the constraints imposed by the tree-child and normality property there are certain restrictions for the
target vertices of the edges we add to the green vertices. We will analyse the parts of the structure. First, since the targets of
the added edges are certainly reticulation vertices, the trees attached to a green vertex cannot be red trees (\emph{cf.} the terminology
at the beginning of Section~\ref{nnt}) and have generating function $\tilde M(z,y)$. We only have to replace $y$ by $y_1$ or $y_2$ or their sum, depending
on whether $g_1$ or $g_2$ or both green vertices, respectively, are allowed to point at the red vertices in this tree (the last situation cannot happen here). The vertices $y$, $g_1$ and $g_2$ cause a
factor $z^3$.

Next, we analyse the contribution of the paths:
\begin{itemize}
\item Path $r$: Both green vertices may point into the attached subtrees, except to their root. The trees are therefore white trees
and the generating function of the path is $1/(1-z\tilde M(z,y_1+y_2)$.
\item Path $k$: The vertex $g_1$ may point to any non-root vertex of the attached trees, $g_2$ may point to any vertex on $k$, except the first one, and any
vertex of the attached trees. Thus the generating function of this path is $P(z,y_2,y_1+y_2,0)$.
\item Path $\ell$: The situation for this path is symmetric to $k$.
\end{itemize}

Overall, this yields the generating function
\begin{align*}
N_{2,2}(z)&=\rcp2\partial_{y_1} \partial_{y_2}\frac{z^3\tilde M(z,y_1)\tilde M(z,y_2)P(z,y_2,y_1+y_2,0)P(z,y_1,y_1+y_2,0)}{1-z\tilde M(z,y_1+y_2)}
\Big\vert_{y_1=0,y_2=0}.
\end{align*}

The exponential generating function for normal networks with
two reticulation vertices is then $N_2(z)=(N_{2,1}(z)+N_{2,2}(z))/4$, where the factor $4$ appears,
because each normal network is generated four times. Simplifying the resulting expression gives \eqref{fctN2}.
\epf

As an easy consequence, we obtain the asymptotic number of networks; see the appendix for numerical data.

\begin{cor}
Let $N_{2,n}$ denote the number of vertex-labeled normal networks with $n$ vertices
and two reticulation vertices. If $n$ is even then $N_{2,n}$ is zero, otherwise
\[
N_{2,n}=n![z^n]N_2(z)=\(\frac{\sqrt2}{e}\)^n n^{n+3}\(\frac{\sqrt{2}}{16}-\frac{3\sqrt{\pi}}{8}\cdot\frac{1}{\sqrt{n}}+{\mathcal O}\left(\frac{1}{n}\right)\),
\]
as $\nti$.
\end{cor}

\bpf This follows by singularity analysis as in the proof of Corollary \ref{asymp-k=1}.\epf

\brem
It turns out that the asymptotic main term is determined by $N_{2,2}(z)$. In hindsight, this is no surprise, because
the corresponding sparsened skeleton has two edges, which leads to three paths made of sequences of trees after all. This
leads to three expressions contributing a singularity in the denominator and together with the number of differentiations this eventually
determines the growth rate of the coefficients of the generating function.
\erem

\subsection{Normal networks with three reticulation vertices}\label{normal-k=3}

In the case of three reticulation vertices we follow the same procedure: We decompose the network
according to how the reticulation vertices are distributed in the network. There are four cases.

\vspace*{0.4cm}
\begin{minipage}{14.2cm}
\begin{enumerate}
\renewcommand{\labelenumi}{Case \arabic{enumi}:}
\item \label{Fall_1} The three reticulation vertices lie on one path, \emph{i.e.}, one reticulation
vertex is ancestor of another, which itself is ancestor of the third one.

\item \label{Fall_2} One reticulation vertex is a common ancestor of the other two, but none of
those two is ancestor of the other one.

\item \label{Fall_3} One reticulation vertex is ancestor of another one, but not of both of them, and the third one is not ancestor of any other reticulation vertex.

\item \label{Fall_4} No reticulation vertex is ancestor of any other reticulation vertex.

\end{enumerate}
\end{minipage}
\vspace*{0.4cm}

\begin{figure}[!h]
	\begin{center}
		\begin{tikzpicture}[
		scale=0.5,
		level/.style={thick},
		virtual/.style={thick,densely dashed},
		trans/.style={thick,<->,shorten >=2pt,shorten <=2pt,>=stealth},
		classical/.style={thin,double,<->,shorten >=4pt,shorten <=4pt,>=stealth}
		]

       \draw (-10.9cm,-1.8cm) node[inner sep=1.5pt,circle,draw,fill] (1) {};
        \draw (-11.3cm,-1.4cm) node[inner sep=1.5pt,circle] (13) {$g_{1}$};
		\draw (-9.9cm,-0.8cm) node[inner sep=1.5pt,circle,draw,fill] (2) {};
		\draw (-8.5cm,0.6cm) node[inner sep=1.5pt,circle,draw,fill] (3) {};
		\draw (-7.5cm,1.6cm) node[inner sep=1.5pt,circle,draw,fill] (4) {};
		\draw (-7.9cm,2cm) node[inner sep=1.5pt,circle] (14) {$g_{2}$};
		\draw (-6.5cm,2.6cm) node[inner sep=1.5pt,circle,draw,fill] (5) {};
		\draw (-5.1cm,4cm) node[inner sep=1.5pt,circle,draw,fill] (6) {};
        \draw (-4.1cm,5cm) node[inner sep=1.5pt,circle,draw,fill] (17) {};
        \draw (-4.5cm,5.4cm) node[inner sep=1.5pt,circle] (22) {$g_3$};
        \draw (-3.1cm,6cm) node[inner sep=1.5pt,circle,draw,fill] (18) {};
        \draw (-1.7cm,7.4cm) node[inner sep=1.5pt,circle,draw,fill] (19) {};
		\draw (-1.3cm,7.9cm) node[inner sep=1.5pt,circle] (15) {$x$};
        \draw (-10.9cm,-3cm) node[regular polygon,regular polygon sides=3,draw,inner sep=0.1cm] (7) {};
		\draw (-9.2cm,-1.5cm) node[regular polygon,regular polygon sides=3,draw,inner sep=0.1cm] (8) {};
		\draw (-7.8cm,-0.1cm) node[regular polygon,regular polygon sides=3,draw,inner sep=0.1cm] (9) {};
		\draw (-5.8cm,1.9cm) node[regular polygon,regular polygon sides=3,draw,inner sep=0.1cm] (11) {};
        \draw (-4.4cm,3.3cm) node[regular polygon,regular polygon sides=3,draw,inner sep=0.1cm] (12) {};
        \draw (-2.4cm,5.3cm) node[regular polygon,regular polygon sides=3,draw,inner sep=0.1cm] (22) {};
        \draw (-1cm,6.7cm) node[regular polygon,regular polygon sides=3,draw,inner sep=0.1cm] (23) {};
        \draw (-9.7cm,-0.6cm) node[inner sep=1.5pt,circle] (13) {};
        \draw (-8.7cm,0.4cm) node[inner sep=1.5pt,circle] (14) {};
        \draw (-6.3cm,2.8cm) node[inner sep=1.5pt,circle] (15) {};
        \draw (-5.3cm,3.8cm) node[inner sep=1.5pt,circle] (16) {};
        \draw (-2.9cm,6.2cm) node[inner sep=1.5pt,circle] (20) {};
        \draw (-1.9cm,7.2cm) node[inner sep=1.5pt,circle] (21) {};
		
		\draw (1)--(2);\draw (3)--(4);\draw (4)--(5);
        \draw (1)--(7);\draw (2)--(8);\draw (3)--(9);
        \draw (5)--(11);\draw (6)--(12);\draw (6)--(17);
        \draw (17)--(18);\draw (18)--(22);\draw (19)--(23);
		\draw[loosely dotted,line width=0.7pt] (13)--(14);
		\draw[loosely dotted,line width=0.7pt] (15)--(16);
        \draw[loosely dotted,line width=0.7pt] (20)--(21);
		\draw[trans] (-10.5cm,-0.6cm) -- (-8.7cm,1.2cm) node[midway,left] {\El{1}};
		\draw[trans] (-7.1cm,2.9cm) -- (-5.3cm,4.6cm) node[midway,left] {\El{2}};
        \draw[trans] (-3.7cm,6.2cm) -- (-1.9cm,8cm) node[midway,left] {\El{3}};

        scale=0.6
        \draw (0cm,0cm) node[inner sep=1.5pt,circle,draw,fill] (1) {};
        \draw (1.2cm,1.2cm) node[inner sep=1.5pt,circle,draw,fill] (2) {};
        \draw (0.8cm,1.6cm) node[inner sep=1.5pt,circle] (26) {$y$};
        \draw (2.4cm,0cm) node[inner sep=1.5pt,circle,draw,fill] (3) {};
        \draw (0.6cm,-0.6cm) node[regular polygon,regular polygon sides=3,draw,inner sep=0.1cm] (4) {};
        \draw (1.8cm,-0.6cm) node[regular polygon,regular polygon sides=3,draw,inner sep=0.1cm] (5) {};
        \draw (-1.6cm,-1.6cm) node[inner sep=1.5pt,circle,draw,fill] (6) {};
        \draw (4cm,-1.6cm) node[inner sep=1.5pt,circle,draw,fill] (7) {};
        \draw (-2.8cm,-2.8cm) node[inner sep=1.5pt,circle,draw,fill] (8) {};
        \draw (-3.2cm,-2.4cm) node[inner sep=1.5pt,circle] (24) {$g_1$};
        \draw (5.2cm,-2.8cm) node[inner sep=1.5pt,circle,draw,fill] (9) {};
        \draw (5.6cm,-2.4cm) node[inner sep=1.5pt,circle] (25) {$g_2$};
        \draw (-1cm,-2.2cm) node[regular polygon,regular polygon sides=3,draw,inner sep=0.1cm] (10) {};
        \draw (3.4cm,-2.2cm) node[regular polygon,regular polygon sides=3,draw,inner sep=0.1cm] (11) {};
        \draw (-2.8cm,-4cm) node[regular polygon,regular polygon sides=3,draw,inner sep=0.1cm] (12) {};
        \draw (5.2cm,-4cm) node[regular polygon,regular polygon sides=3,draw,inner sep=0.1cm] (13) {};
        \draw (2.4cm,2.4cm) node[inner sep=1.5pt,circle,draw,fill] (14) {};
        \draw (4cm,4cm) node[inner sep=1.5pt,circle,draw,fill] (15) {};
        \draw (5.2cm,5.2cm) node[inner sep=1.5pt,circle,draw,fill] (28) {};
        \draw (4.8cm,5.6cm) node[inner sep=1.5pt,circle] (27) {$g_3$};
        \draw (6.4cm,6.4cm) node[inner sep=1.5pt,circle,draw,fill] (29) {};
        \draw (8cm,8cm) node[inner sep=1.5pt,circle,draw,fill] (30) {};
        \draw (8.4cm,8.4cm) node[inner sep=1.5pt,circle] (31) {$x$};
        \draw (3cm,1.8cm) node[regular polygon,regular polygon sides=3,draw,inner sep=0.1cm] (16) {};
        \draw (4.6cm,3.4cm) node[regular polygon,regular polygon sides=3,draw,inner sep=0.1cm] (17) {};
        \draw (7cm,5.8cm) node[regular polygon,regular polygon sides=3,draw,inner sep=0.1cm] (34) {};
        \draw (8.6cm,7.4cm) node[regular polygon,regular polygon sides=3,draw,inner sep=0.1cm] (35) {};
        \draw (-1.4cm,-1.4cm) node[inner sep=1.5pt,circle] (18) {};
        \draw (-0.2cm,-0.2cm) node[inner sep=1.5pt,circle] (19) {};
        \draw (3.8cm,-1.4cm) node[inner sep=1.5pt,circle] (20) {};
        \draw (2.6cm,-0.2cm) node[inner sep=1.5pt,circle] (21) {};
        \draw (2.6cm,2.6cm) node[inner sep=1.5pt,circle] (22) {};
        \draw (3.8cm,3.8cm) node[inner sep=1.5pt,circle] (23) {};
        \draw (6.6cm,6.6cm) node[inner sep=1.5pt,circle] (32) {};
        \draw (7.8cm,7.8cm) node[inner sep=1.5pt,circle] (33) {};

        \draw (8)--(12);\draw (8)--(6);\draw (6)--(10);\draw (1)--(4);\draw (1)--(2);
        \draw (2)--(3);\draw (3)--(5);\draw (7)--(11);\draw (7)--(9);\draw (9)--(13);
        \draw (2)--(14);\draw (14)--(16);\draw (15)--(17);\draw (15)--(28);\draw (28)--(29);
        \draw (29)--(34);\draw (30)--(35);
        \draw [loosely dotted,line width=0.7pt] (18)--(19);
        \draw [loosely dotted,line width=0.7pt] (20)--(21);
        \draw [loosely dotted,line width=0.7pt] (22)--(23);
        \draw [loosely dotted,line width=0.7pt] (32)--(33);
        \draw[trans] (-2.2cm,-1.4cm) -- (-0.2cm,0.6cm) node[midway,left] {\El{1}};
		\draw[trans] (2.6cm,0.6cm) -- (4.6cm,-1.4cm) node[midway,right] {\El{2}};
        \draw[trans] (1.8cm,2.6cm) -- (3.8cm,4.6cm) node[midway,left] {\El{3}};
        \draw[trans] (5.8cm,6.6cm) -- (7.8cm,8.6cm) node[midway,left] {\El{4}};
\end{tikzpicture}
	\end{center}
\caption{Two of the four possible structures of Motzkin skeletons of networks with three reticulation vertices. The left one arises from the sparsened skeleton which is a path of length $2$ and the right one arises from a unary vertex to which a cherry is attached.}
\label{fig_skeleton_3}
	\end{figure}
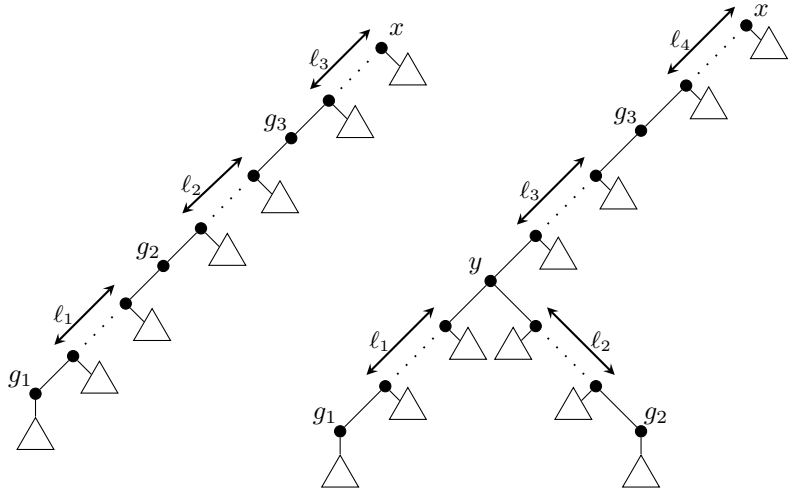

\begin{figure}
		\begin{center}
			\begin{tikzpicture}[
			scale=0.5,
			level/.style={thick},
			virtual/.style={thick,densely dashed},
			trans/.style={thick,<->,shorten >=2pt,shorten <=2pt,>=stealth},
			classical/.style={thin,double,<->,shorten >=4pt,shorten <=4pt,>=stealth}
			]
			
	    \draw (0cm,0cm) node[inner sep=1.5pt,circle,draw,fill] (1) {};
        \draw (1.2cm,1.2cm) node[inner sep=1.5pt,circle,draw,fill] (2) {};
        \draw (0.8cm,1.6cm) node[inner sep=1.5pt,circle] (26) {$z$};
        \draw (2.4cm,0cm) node[inner sep=1.5pt,circle,draw,fill] (3) {};
        \draw (0.6cm,-0.6cm) node[regular polygon,regular polygon sides=3,draw,inner sep=0.1cm] (4) {};
        \draw (1.8cm,-0.6cm) node[regular polygon,regular polygon sides=3,draw,inner sep=0.1cm] (5) {};
        \draw (-1.6cm,-1.6cm) node[inner sep=1.5pt,circle,draw,fill] (6) {};
        \draw (4cm,-1.6cm) node[inner sep=1.5pt,circle,draw,fill] (7) {};
        \draw (-2.8cm,-2.8cm) node[inner sep=1.5pt,circle,draw,fill] (8) {};
        \draw (-3.2cm,-2.4cm) node[inner sep=1.5pt,circle] (24) {$g_1$};
        \draw (5.2cm,-2.8cm) node[inner sep=1.5pt,circle,draw,fill] (9) {};
        \draw (5.6cm,-2.4cm) node[inner sep=1.5pt,circle] (25) {$g_2$};
        \draw (-1cm,-2.2cm) node[regular polygon,regular polygon sides=3,draw,inner sep=0.1cm] (10) {};
        \draw (3.4cm,-2.2cm) node[regular polygon,regular polygon sides=3,draw,inner sep=0.1cm] (11) {};
        \draw (-2.8cm,-4cm) node[regular polygon,regular polygon sides=3,draw,inner sep=0.1cm] (12) {};
        \draw (5.2cm,-4cm) node[regular polygon,regular polygon sides=3,draw,inner sep=0.1cm] (13) {};
        \draw (2.4cm,2.4cm) node[inner sep=1.5pt,circle,draw,fill] (14) {};
        \draw (4cm,4cm) node[inner sep=1.5pt,circle,draw,fill] (15) {};
        \draw (5.2cm,5.2cm) node[inner sep=1.5pt,circle,draw,fill] (28) {};
        \draw (4.8cm,5.6cm) node[inner sep=1.5pt,circle] (42) {$y$};
        \draw (6.4cm,6.4cm) node[inner sep=1.5pt,circle,draw,fill] (29) {};
        \draw (8cm,8cm) node[inner sep=1.5pt,circle,draw,fill] (30) {};
        \draw (8.4cm,8.4cm) node[inner sep=1.5pt,circle] (27) {$x$};
        \draw (3cm,1.8cm) node[regular polygon,regular polygon sides=3,draw,inner sep=0.1cm] (16) {};
        \draw (4.6cm,3.4cm) node[regular polygon,regular polygon sides=3,draw,inner sep=0.1cm] (17) {};
        \draw (6.4cm,4cm) node[inner sep=1.5pt,circle,draw,fill] (33) {};
        \draw (8cm,2.4cm) node[inner sep=1.5pt,circle,draw,fill] (34) {};
        \draw (9.2cm,1.2cm) node[inner sep=1.5pt,circle,draw,fill] (35) {};
        \draw (9.6cm,1.6cm) node[inner sep=1.5pt,circle] (36) {$g_3$};
        \draw (5.8cm,3.4cm) node[regular polygon,regular polygon sides=3,draw,inner sep=0.1cm] (41) {};
        \draw (7.4cm,1.8cm) node[regular polygon,regular polygon sides=3,draw,inner sep=0.1cm] (40) {};
        \draw (9.2cm,0cm) node[regular polygon,regular polygon sides=3,draw,inner sep=0.1cm] (39) {};
        \draw (-1.4cm,-1.4cm) node[inner sep=1.5pt,circle] (18) {};
        \draw (-0.2cm,-0.2cm) node[inner sep=1.5pt,circle] (19) {};
        \draw (3.8cm,-1.4cm) node[inner sep=1.5pt,circle] (20) {};
        \draw (2.6cm,-0.2cm) node[inner sep=1.5pt,circle] (21) {};
        \draw (2.6cm,2.6cm) node[inner sep=1.5pt,circle] (22) {};
        \draw (3.8cm,3.8cm) node[inner sep=1.5pt,circle] (23) {};
        \draw (6.6cm,6.6cm) node[inner sep=1.5pt,circle] (31) {};
        \draw (7.8cm,7.8cm) node[inner sep=1.5pt,circle] (32) {};
        \draw (6.6cm,3.8cm) node[inner sep=1.5pt,circle] (37) {};
        \draw (7.8cm,2.6cm) node[inner sep=1.5pt,circle] (38) {};
        \draw (7cm,5.8cm) node[regular polygon,regular polygon sides=3,draw,inner sep=0.1cm] (43) {};
        \draw (8.6cm,7.4cm) node[regular polygon,regular polygon sides=3,draw,inner sep=0.1cm] (44) {};

        \draw (8)--(12);\draw (8)--(6);\draw (6)--(10);\draw (1)--(4);\draw (1)--(2);
        \draw (2)--(3);\draw (3)--(5);\draw (7)--(11);\draw (7)--(9);\draw (9)--(13);
        \draw (2)--(14);\draw (14)--(16);\draw (15)--(17);\draw (15)--(28);\draw (28)--(29);
        \draw (28)--(33);\draw (34)--(35);\draw (33)--(41);\draw (34)--(40);\draw (35)--(39);
        \draw (29)--(43);\draw (30)--(44);
        \draw [loosely dotted,line width=0.7pt] (18)--(19);
        \draw [loosely dotted,line width=0.7pt] (20)--(21);
        \draw [loosely dotted,line width=0.7pt] (22)--(23);
        \draw [loosely dotted,line width=0.7pt] (31)--(32);
        \draw [loosely dotted,line width=0.7pt] (37)--(38);
        \draw[trans] (-2.2cm,-1.4cm) -- (-0.2cm,0.6cm) node[midway,left] {\El{1}};
		\draw[trans] (2.6cm,0.6cm) -- (4.6cm,-1.4cm) node[midway,right] {\El{2}};
        \draw[trans] (1.8cm,2.6cm) -- (3.8cm,4.6cm) node[midway,left] {\El{4}};
        \draw[trans] (5.8cm,6.6cm) -- (7.8cm,8.6cm) node[midway,left] {\El{5}};
        \draw[trans] (6.6cm,4.6cm) -- (8.6cm,2.6cm) node[midway,right] {\El{3}};
			\end{tikzpicture}
		\end{center}
		\caption{The fourth possible structure of Motzkin skeletons of
networks with three reticulation vertices. It arises from a sparsened skeleton which is a rooted binary caterpillar with three leaves.}
\label{fig_skeleton_5}
	\end{figure}
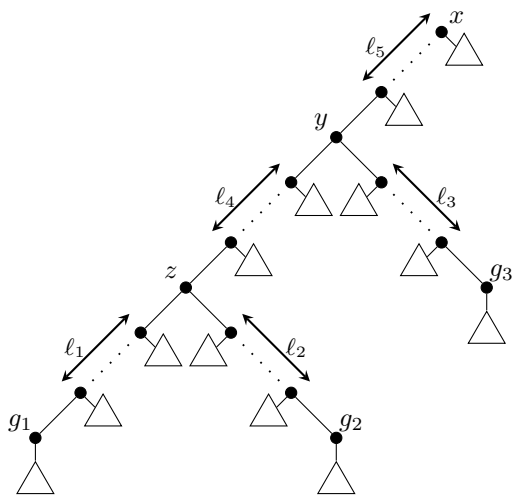

The Motzkin skeletons of the first two cases are depicted in
Figure~\ref{fig_skeleton_3}, Case 3 is depicted in Figure~\ref{fig_skeleton_4} and the last case is depicted
in Figure~\ref{fig_skeleton_5}. For all those cases, the
lengths of the paths connecting two reticulation vertices or connecting a reticulation vertex with the
last common ancestor of two reticulation vertices are the free parameters. To each vertex of such a
path we may attach a Motzkin tree which must be shaped in such a way that the condition for
normality of the network is respected. So, we will set up generating functions $N_{3,1}(z)$,
$N_{3,2}(z)$, $N_{3,3}(z)$, $N_{3,4}(z)$ associated to Motzkin skeletons of the four respective
cases and since the procedure will generate each normal network eight times, the generating
function of normal networks with three reticulation vertices is
\[
N_3(z)=\frac{N_{3,1}(z)+N_{3,2}(z)+N_{3,3}(z)+N_{3,4}(z)}8.
\]

We start with Case~\ref{Fall_1}, see Figure~\ref{fig_skeleton_3}, left tree.
As in the previous section we call the path from the root to the parent of $g_1$ the spine of the
Motzkin skeleton. Then $g_1$ can point to each non-root vertex of the subtrees attached to any of the
vertices of the spine. Likewise, the pointing options for $g_2$ are the non-root vertices of the
subtrees attached to the vertices of the sub-paths $\ell_2\cup \ell_3$ of the spine. The situation for $g_3$ is analogous.

From this we obtain the following expression for the generating function of all normal networks
with three reticulation vertices which are distributed within the network as shown in
Figure~\ref{fig_skeleton_3}, left:
\begin{align*}
N_{3,1}(z)&=\Y\left(
\dfrac{z^3M(z,0)}{(1-z\tilde{M}(z,y_{1}))
(1-z\tilde{M}(z,y_{1}+y_{2}))(1-z\tilde{M}(z,y_{1}+y_{2}+y_{3}))}
\right),
\end{align*}
where $\tilde{M}(z,y)$ is as in the last subsection (\emph{cf.}~\eqref{tree-without-root}) and $\Y$ denotes the operator differentiating with respect to $y_1,\ y_2,\ y_3$ and setting
$y_1=y_2=y_3=0$ afterwards, \emph{i.e.}, $\Y
f(z,y_1,y_2,y_3)=\(\partial_{y_1}\partial_{y_2} \partial_{y_3}f\)(z,0,0,0)$.

Next we will determine the generating function of all normal networks belonging to Case~\ref{Fall_2},
which have Motzkin skeletons as shown on the right of Figure~\ref{fig_skeleton_3}. As in the previous section we analyse the substructures. There are four
vertices in the sparsened skeleton, yielding a factor $z^4$. The red vertices in the (white) subtree attached to $g_1$ may only be targets of the edge coming from
$g_2$, for the subtree attached to $g_2$ \emph{vice versa}.

\begin{itemize}
\item Paths $\ell_3$ and $\ell_4$: These paths are sequences of vertices, each with a white subtree attached to it. For $\ell_4$ each green vertex is allowed to point at the red vertices in these white subtrees. Pointing to the vertices of the path is not allowed. Likewise, the corresponding vertices in the subtrees of
$\ell_3$ are forbidden for $g_3$ by the normality condition.
\item Paths $\ell_1$ and $\ell_2$: They are symmetric, so we discuss $\ell_1$. The non-root vertices of the subtrees are the only allowed targets
for the edge from $g_1$. The edge from $g_2$ may end at each vertex of the subtrees and the vertices of the path itself except the first vertex. There are no options
for $g_3$. So, we must distinguish three classes: the first vertex of the path, the other vertices of the path and the roots of the trees, the vertices ``strictly inside'' the trees. These precisely correspond to the variables $\hat y$, $y$ and $\tilde y$, respectively, of the function $P(z,y,\tilde y,\hat y)$ introduced in \eqref{PGF}. We obtain $P(z,y_2,y_1+y_2,0)$ for $\ell_1$.
\end{itemize}

Overall, this gives, again using the operator $\Y$ defined above, the generating function 	
	
\begin{align*}
N_{3,2}(z)&=\rcp2\Y\left(\frac{z^4\tilde M(z,y_1)\tilde M(z,y_2)P(z,y_1,y_1+y_2,0)P(z,y_2,y_1+y_2,0)}{(1-z\tilde M(z,y_1+y_2+y_3))(1-z\tilde M(z,y_1+y_2))}\right)
\end{align*}	

Case~\ref{Fall_3} is the one shown in Figure~\ref{fig_skeleton_4}. The sparsened skeleton has 4 vertices and the subtrees attached to $g_1$ and $g_3$ are
white trees. The red vertices of the subtree of $g_1$ may be targeted by the edges starting either in $g_2$ or in $g_3$, the red vertices of the other tree by
edges from $g_1$.

Next we inspect the paths:
\begin{itemize}
\item Path $\ell_4$: All green vertices may point to the non-root vertices of the (white) subtrees. Pointing to the path itself is not allowed.
\item Path $\ell_3$: The edge starting at $g_3$ may point to non-root vertices of the subtrees, but neither to the roots nor to the vertices of the path itself.
There is no option for $g_2$. All but the first vertex of the path as well as all tree vertices can be the end point of the edge starting at $g_1$.
\item Path $\ell_1$: Similar to $\ell_3$. The edges from $g_2$ and $g_3$ may point anywhere except to the first vertex of the path. The non-root vertices of the
subtrees may be targeted by $g_1$ as well.
\item Path $\ell_2$: All green vertices may point to the non-root vertices of the subtrees. To point at the vertices on the path or to the root vertices of the subtrees is only allowed for $g_1$. Again, the first vertex of the path is the exception. It must not be red by the normality condition.
\end{itemize}

Altogether, we obtain for the generating function $N_{3,3}(z)$ of Case~\ref{Fall_3} the expression
\begin{align*}
	N_{3,3}(z)&=
	\Y\left(\frac{z^4\tilde M(z,y_2+y_3)\tilde M(z,y_1)}{1-z\tilde M(z,y_1+y_2+y_3)} P(z,y_1,y_1+y_3,0) \right. \\
	&\ \ \qquad \times P(z,y_2+y_3,y_1+y_2+y_3,0) P(z,y_1,y_1+y_2+y_3,0)\Bigg).
\end{align*}

\begin{figure}
	\begin{center}
		\begin{tikzpicture}[
		scale=0.5,
		level/.style={thick},
		virtual/.style={thick,densely dashed},
		trans/.style={thick,<->,shorten >=2pt,shorten <=2pt,>=stealth},
		classical/.style={thin,double,<->,shorten >=4pt,shorten <=4pt,>=stealth}
		]

        \draw (0cm,0cm) node[inner sep=1.5pt,circle,draw,fill] (1) {};
        \draw (1.2cm,1.2cm) node[inner sep=1.5pt,circle,draw,fill] (2) {};
        \draw (0.8cm,1.6cm) node[inner sep=1.5pt,circle] (26) {$y$};
        \draw (2.4cm,0cm) node[inner sep=1.5pt,circle,draw,fill] (3) {};
        \draw (0.6cm,-0.6cm) node[regular polygon,regular polygon sides=3,draw,inner sep=0.1cm] (4) {};
        \draw (1.8cm,-0.6cm) node[regular polygon,regular polygon sides=3,draw,inner sep=0.1cm] (5) {};
        \draw (-1.6cm,-1.6cm) node[inner sep=1.5pt,circle,draw,fill] (6) {};
        \draw (4cm,-1.6cm) node[inner sep=1.5pt,circle,draw,fill] (7) {};
        \draw (-2.8cm,-2.8cm) node[inner sep=1.5pt,circle,draw,fill] (8) {};
        \draw (-3.2cm,-2.4cm) node[inner sep=1.5pt,circle] (24) {$g_1$};
        \draw (5.2cm,-2.8cm) node[inner sep=1.5pt,circle,draw,fill] (9) {};
        \draw (5.6cm,-2.4cm) node[inner sep=1.5pt,circle] (25) {$g_2$};
        \draw (6.4cm,-4cm) node[inner sep=1.5pt,circle,draw,fill] (28) {};
        \draw (8cm,-5.6cm) node[inner sep=1.5pt,circle,draw,fill] (31) {};
        \draw (9.2cm,-6.8cm) node[inner sep=1.5pt,circle,draw,fill] (29) {};
        \draw (9.6cm,-6.4cm) node[inner sep=1.5pt,circle] (35) {$g_3$};
        \draw (-1cm,-2.2cm) node[regular polygon,regular polygon sides=3,draw,inner sep=0.1cm] (10) {};
        \draw (3.4cm,-2.2cm) node[regular polygon,regular polygon sides=3,draw,inner sep=0.1cm] (11) {};
        \draw (-2.8cm,-4cm) node[regular polygon,regular polygon sides=3,draw,inner sep=0.1cm] (12) {};
        \draw (9.2cm,-8cm) node[regular polygon,regular polygon sides=3,draw,inner sep=0.1cm] (13) {};
        \draw (5.8cm,-4.6cm) node[regular polygon,regular polygon sides=3,draw,inner sep=0.1cm] (30) {};
        \draw (7.4cm,-6.2cm) node[regular polygon,regular polygon sides=3,draw,inner sep=0.1cm] (32) {};
        \draw (2.4cm,2.4cm) node[inner sep=1.5pt,circle,draw,fill] (14) {};
        \draw (4cm,4cm) node[inner sep=1.5pt,circle,draw,fill] (15) {};
        \draw (4.4cm,4.4cm) node[inner sep=1.5pt,circle] (27) {$x$};
        \draw (3cm,1.8cm) node[regular polygon,regular polygon sides=3,draw,inner sep=0.1cm] (16) {};
        \draw (4.6cm,3.4cm) node[regular polygon,regular polygon sides=3,draw,inner sep=0.1cm] (17) {};
        \draw (-1.4cm,-1.4cm) node[inner sep=1.5pt,circle] (18) {};
        \draw (-0.2cm,-0.2cm) node[inner sep=1.5pt,circle] (19) {};
        \draw (3.8cm,-1.4cm) node[inner sep=1.5pt,circle] (20) {};
        \draw (2.6cm,-0.2cm) node[inner sep=1.5pt,circle] (21) {};
        \draw (2.6cm,2.6cm) node[inner sep=1.5pt,circle] (22) {};
        \draw (3.8cm,3.8cm) node[inner sep=1.5pt,circle] (23) {};
        \draw (6.6cm,-4.2cm) node[inner sep=1.5pt,circle] (33) {};
        \draw (7.8cm,-5.4cm) node[inner sep=1.5pt,circle] (34) {};

        \draw (8)--(12);\draw (8)--(6);\draw (6)--(10);\draw (1)--(4);\draw (1)--(2);
        \draw (2)--(3);\draw (3)--(5);\draw (7)--(11);\draw (7)--(9);\draw (29)--(13);
        \draw (2)--(14);\draw (14)--(16);\draw (15)--(17);\draw (9)--(28);\draw (28)--(30);
        \draw (31)--(32);\draw (31)--(29);
        \draw [loosely dotted,line width=0.7pt] (18)--(19);
        \draw [loosely dotted,line width=0.7pt] (20)--(21);
        \draw [loosely dotted,line width=0.7pt] (22)--(23);
        \draw [loosely dotted,line width=0.7pt] (33)--(34);
        \draw[trans] (-2.2cm,-1.4cm) -- (-0.2cm,0.6cm) node[midway,left] {\El{1}};
		\draw[trans] (2.6cm,0.6cm) -- (4.6cm,-1.4cm) node[midway,right] {\El{2}};
        \draw[trans] (1.8cm,2.6cm) -- (3.8cm,4.6cm) node[midway,left] {\El{4}};
        \draw[trans] (6.6cm,-3.4cm) -- (8.6cm,-5.4cm) node[midway,right] {\El{3}};
\end{tikzpicture}
	\end{center}
\caption{The third possible structure of Motzkin skeletons of networks with three reticulation vertices. It arises from the sparsened skeleton which consists of a root with a left child and path of length $2$ as right subtree.}
\label{fig_skeleton_4}
	\end{figure}
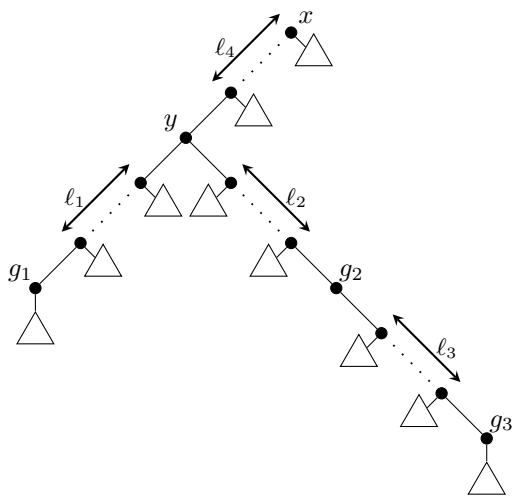

The last case of normal networks has Motzkin skeletons as shown in Figure~\ref{fig_skeleton_5}.
The restriction for the target vertex of the edges to be added at $g_1$, $g_2$ and $g_3$ follow the
analogous rules in order to meet the normality constraint. Setting up the generating function follows the same pattern as before.
We omit now the details and get from path analysis after all

\begin{align*}
&N_{3,4}(z)=
\\
&\rcp2 \Y\left(
\frac{z^5\tilde M(z,y_1+y_2)\tilde M(z,y_1+y_3)\tilde M(z,y_2+y_3)}{1-z\tilde M(z,y_1+y_2+y_3)}P(z,y_1+y_2,y_1+y_2+y_3,0)
\right.
\\
&\qquad\times  P(z,y_1+y_3,y_1+y_2+y_3,y_3)P(z,y_2+y_3,y_1+y_2+y_3,y_3)P(z,y_3,y_1+y_2+y_3,0)\Bigg).
\end{align*}

Overall, by collecting everything, we obtain the following result.

\bp\label{exp-k=3}
The exponential generating function for vertex-labeled normal networks with
three reticulation vertices is
\begin{equation*}
N_3(z)=z\frac{\tilde{a}_3(z^2)-\tilde{b}_3(z^2)\sqrt{1-2z^2}}{(1-2z^2)^{11/2}},
\end{equation*}
where
\begin{align*}
\tilde{a}_3(z)&=270z^6-\frac{2187}{2}z^5+576 z^4-9 z^3
\end{align*}
and
\begin{align*}
\tilde{b}_3(z)&=18z^6-531z^5+567z^4-9z^3
\end{align*}
\ep

As a consequence we obtain the following result; see the appendix for numerical data.
\begin{cor}
Let $N_{3,n}$ denote the number of vertex-labeled normal networks with $n$ vertices
and three reticulation vertices. If $n$ is even then $N_{3,n}$ is zero, otherwise
\[
N_{3,n}=n![z^n]N_3(z)=\(\frac{\sqrt2}{e}\)^n n^{n+5}\(\frac{\sqrt{2}}{192}-\frac{3\sqrt{\pi}}{64}\cdot\frac{1}{\sqrt{n}}+{\mathcal O}\left(\frac{1}{n}\right)\),
\]
as $\nti$.
\end{cor}

\bpf This follows by singularity analysis as for $k=1$ and $k=2$ above.\epf

\subsection{Normal networks with a fixed number of reticulation vertices}

By looking at Proposition \ref{exp-k=1}, Proposition \ref{exp-k=2} and Proposition \ref{exp-k=3}, one clearly sees a pattern for the exponential generating function of normal networks. In this section, we will prove that this pattern continues to hold for the exponential generating function of normal networks with $k$ reticulation vertices. This will then be used to prove the remaining claims of Theorem \ref{main-thm-1}.

We start with a technical lemma. Therefore, consider the following function
\begin{equation}\label{func-G}
G(z,y)=\frac{a(z,y)-b(z,y)\sqrt{1-2z^2-4yz^3}}{1+2yz},
\end{equation}
where $a(z,y),b(z,y)$ are polynomials in $z$ and $y$ with $a(z,0)=b(z,0)=1$. This function will be used as a building block for construction the exponential generating
function of normal networks. We need the following simple properties of this function.
\bl\label{tech-lmm}
\begin{itemize}
\item[(a)] For all $\ell\geq 1$,
\[
\frac{\partial^{\ell}}{\partial
y^{\ell}}G(z,y)\Big\vert_{y=0}=\frac{c_{\ell}(z)-d_{\ell}(z)\sqrt{1-2z^2}}{(1-2z^2)^{\ell-1/2}},
\]
where $c_{\ell}(z)$ and $d_{\ell}(z)$ are suitable polynomials.
\item[(b)] For all $\ell\geq 0$,
\[
\frac{\partial^{\ell}}{\partial
y^{\ell}}\frac{1}{1-G(z,y)}
\Big\vert_{y=0}=\frac{e_{\ell}(z)-f_{\ell}(z)\sqrt{1-2z^2}}{(1-2z^2)^{\ell+1/2}},
\]
where $e_{\ell}(z)$ and $f_{\ell}(z)$ are suitable polynomials.
\end{itemize}
\el
\bpf For the proof of part (a), by differentiation
\[
\frac{\partial^{\ell}}{\partial y^{\ell}}G(z,y)=\frac{a_{\ell}(z,y)-b_{\ell}(z,y)\sqrt{1-2z^2-4yz^3}}{(1+2yz)^{\ell+1}(1-2z^2-4yz^3)^{\ell-1/2}}
\]
with suitable polynomials $a_{\ell}(z,y)$ an $b_{\ell}(z,y)$. (Note that this becomes incorrect for $\ell=0$). The claim follows now by setting $y=0$.

For the proof of part (b), we use induction. Note that $\ell=0$ is trivial. Now, assume that the claim holds for all $\tilde{\ell}<\ell$. Then, by Leibnitz rule
\begin{align*}
&\frac{\partial^{\ell}}{\partial y^{\ell}}\frac{1}{1-G(z,y)}\Big\vert_{y=0}=\frac{\partial^{\ell-1}}{\partial y^{\ell-1}}\left(\frac{1}{1-G(z,y)}\cdot\frac{1}{1-G(z,y)}\cdot G'(z,y)\right)\Big\vert_{y=0}\\
&=\sum_{k_1+k_2+k_3=\ell-1}\binom{\ell-1}{k_1,k_2,k_3}\frac{\partial^{k_1}}{\partial y^{k_1}}\frac{1}{1-G(z,y)}\Big\vert_{y=0}\cdot\frac{\partial^{k_2}}{\partial y^{k_2}}\frac{1}{1-G(z,y)}\Big\vert_{y=0}\cdot G^{(k_3+1)}(z,y)\Big\vert_{y=0}.
\end{align*}
Plugging into this the induction hypothesis and part (a) gives the claimed form with power of the denominator equal to
\[
k_1+1/2+k_2+1/2+k_3+1/2=\ell+1/2.
\]
This proves the result.
\epf

Now, we can prove the following result which generalizes Proposition \ref{exp-k=1}, Proposition \ref{exp-k=2} and Proposition \ref{exp-k=3}.

\bp\label{gen-normal-k}
The exponential generating function for vertex-labeled normal networks with $k$ reticulation vertices is
\[
N_{k}(z)=\frac{a_k(z)-b_k(z)\sqrt{1-2z^2}}{(1-2z^2)^{2k-1/2}},
\]
where $a_k(z)$ and $b_k(z)$ are suitable polynomials.
\ep
\bpf Fix a type of Motzkin skeletons (arising from a sparsened skeleton) for generating normal networks with $k$ reticulation vertices. As explained in the cases $k=1,2,3$, the exponential generating function of the normal networks arising from these skeletons is a product of generating functions for the paths which are either counted by $1/(1-z\tilde{M})$ or $P$ multiplied with a $z$ for each vertex of the sparsened skeleton and the generating functions of the Motzkin trees attached to the leaves. In particular note that $z\tilde{M}$ is of the form (\ref{func-G}) and the denominators of $P$ is one minus a function of the form (\ref{func-G}). Also, note that all these functions $G$ have polynomials satisfying $a(z,0)=b(z,0)=1$.

In summary, we have that the exponential generating function $N_{k}(z)$ for normal networks is a sum of terms of the form
\begin{equation}\label{terms-N}
\partial_{y_1}\cdots\partial_{y_k}\frac{G_1(z,y)\cdots G_{s}(z,y)}{(1-G_{s+1}(z,y))\cdots(1-G_{s+t}(z,y))}\Big\vert_{y_1=0,\ldots,y_k=0},
\end{equation}
where the number of functions $G_{s+i}(z,y)$ is bounded by the number of edges of the sparsened skeleton increased by one (for the sequence of trees added above the root when constructing the Motzkin skeletons). Moreover, $y$ is the sum of the $y_i$'s where not all of them must be present and the missing ones can differ from one occurrence to the next in the above formula. Note that because of this special form of $y$, we can apply the above lemma after expanding (\ref{terms-N}) and obtain that
\begin{equation}\label{gen-form}
N_k(z)=\frac{a_k(z)-b_k(z)\sqrt{1-2z^2}}{(1-2z^2)^{p}}.
\end{equation}

What remains is to show that $p=2k-1/2$. For this observe that (\ref{terms-N}) without the
derivatives is of the general form given in (\ref{gen-form}) with the exponent of the denominator
equals $t/2$ which reaches its maximum for the sparsened skeleton with the maximal number of edges and is thus at most $k-1/2$. Also, from the above lemma, we see that each differentiation increases the exponent
by $1$. Thus, the exponent of (\ref{terms-N}) when written as (\ref{gen-form}) is at most $2k-1/2$. Adding up this terms gives the claim.
\epf

\begin{cor}\label{cor-gen}
We have
\[
N_k(z)=z\frac{\tilde{a}_k(z^2)-\tilde{b}_k(z^2)\sqrt{1-2z^2}}{(1-2z^2)^{2k-1/2}},
\]
where $\tilde{a}_k(z)$ and $\tilde{b}_k(z)$ are suitable polynomials.
\end{cor}
\bpf Observe that $N_k(-z)=-N_k(z)$ since phylogenetic networks necessarily have an odd number of vertices.
Thus,
\[
\frac{a_k(-z)-b_k(-z)\sqrt{1-2z^2}}{(1-2z^2)^{2k-1/2}}=
-\frac{a_k(z)-b_k(z)\sqrt{1-2z^2}}{(1-2z^2)^{2k-1/2}}.
\]
This implies
\[
a_k(-z)+a_k(z)=(b_k(-z)+b_k(z))\sqrt{1-2z^2}
\]
which is obviously only possible if
\[
a_k(-z)=-a_k(z)\qquad\text{and}\qquad b_k(-z)=-b_k(z),
\]
\emph{i.e.}, both are odd functions. From this the result follows.
\epf

Now, we can finish the proof of Theorem \ref{main-thm-1}.
\begin{cor}\label{main-thm}
Let $N_{k,n}$ denote the number of vertex-labeled normal networks with $n$ vertices and $k$ reticulation vertices. If $n$ is even then $N_{k,n}$ is zero, otherwise there is a positive constant $\tilde{c}_{k}$ such that
\[
N_{k,n}=n![z^n]N_k(z)\sim \tilde{c}_k\left(\frac{\sqrt{2}}{e}\right)^n n^{n+2k-1},
\]
$n\rightarrow\infty$.
\end{cor}
\bpf From the above corollary,
\[
N_{n,k}=n![z^n]z\frac{\tilde{a}_k(z^2)-\tilde{b}_k(z^2)\sqrt{1-2z^2}}{(1-2z^2)^{2k-1/2}}.
\]
From this, by singularity analysis and Stirling's formula, the claimed expansion follows with
\[
\tilde{c}_k=\frac{2\sqrt{2\pi}\tilde{a}_k(1/2)}{4^k\Gamma(2k-1/2)}.
\]

What is left is to prove that $\tilde{c}_k>0$ (note that we already showed this for $k=1,2,3$ directly).
This will follow from Proposition~\ref{prop_5} below which shows that already a subset of the set of normal
networks with $k$ reticulation vertices satisfies the above claimed asymptotics with a
positive constant.
\epf

The proof of Corollary~\ref{main-thm} relies on the fact that a certain constant (called $\tilde
c_k$ there) is positive. This constant is related to the number of normal phylogenetic networks;
it is the multiplicative constant of the asymptotic main term. We will construct a subclass of the
class of normal networks and show that the number of networks in that subclass is the same as for
normal networks up to a positive multiplicative constant. The result is presented in
Proposition~\ref{prop_5} below and closes the small gap left in the proof of
Corollary~\ref{main-thm}.

For this purpose, we consider all the normal networks which are generated (possibly with
duplicity) from a sparsened skeleton which is a rooted binary caterpillar, \emph{i.e.}, a sparsened skeleton of the form
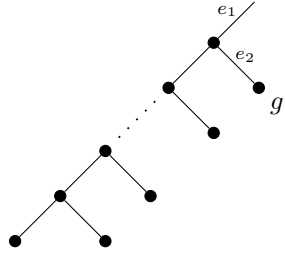
\begin{center}
\begin{tikzpicture}[scale=0.6,
		level/.style={thick},
		virtual/.style={thick,densely dashed},
		trans/.style={thick,<->,shorten >=2pt,shorten <=2pt,>=stealth},
		classical/.style={thin,double,<->,shorten >=4pt,shorten <=4pt,>=stealth}]

\draw (0cm,0cm) node[inner sep=1.5pt,circle,draw,fill] (1) {};
\draw (1cm,1cm) node[inner sep=1.5pt,circle,draw,fill] (2) {};
\draw (2cm,2cm) node[inner sep=1.5pt,circle,draw,fill] (3) {};
\draw (2cm,0cm) node[inner sep=1.5pt,circle,draw,fill] (4) {};
\draw (3cm,1cm) node[inner sep=1.5pt,circle,draw,fill] (5) {};
\draw (3.4cm,3.4cm) node[inner sep=1.5pt,circle,draw,fill] (6) {};
\draw (4.4cm,4.4cm) node[inner sep=1.5pt,circle,draw,fill] (7) {};
\draw (5.4cm,5.4cm) node[inner sep=1.5pt,circle] (8) {};
\draw (4.4cm,2.4cm) node[inner sep=1.5pt,circle,draw,fill] (9) {};
\draw (5.4cm,3.4cm) node[inner sep=1.5pt,circle,draw,fill] (10) {};
\draw (4.7cm,5.1cm) node[inner sep=1.5pt,circle] (11) {$\scriptstyle{e_1}$};
\draw (5.1cm,4.1cm) node[inner sep=1.5pt,circle] (12) {$\scriptstyle{e_2}$};
\draw (5.8cm,3cm) node[inner sep=1.5pt,circle] (13) {$g$};
\draw (2.2cm,2.2cm) node[inner sep=1.5pt,circle] (11) {};
\draw (3.2cm,3.2cm) node[inner sep=1.5pt,circle] (12) {};

\draw (1)--(2);\draw (2)--(3);\draw (2)--(4);\draw (3)--(5);
\draw (6)--(7);\draw (7)--(8);\draw (6)--(9);\draw (7)--(10);
\draw[loosely dotted,line width=0.7pt] (11)--(12);
\end{tikzpicture}
\end{center}

\vspace*{0.2cm}
\noindent (For the discussion below, we have added an edge from the root.) Note that by the same arguments as above, these networks are also counted by an exponential generating function of the form
\begin{equation}\label{gen-func-cat}
C_k(z)=z\frac{\tilde{e}_k(z^2)-\tilde{f}_k(z^2)\sqrt{1-2z^2}}{(1-2z^2)^{2k-1/2}},
\end{equation}
where $\tilde{e}_k(z)$ and $\tilde{f}_k(z)$ are suitable polynomials.

Now, we are in position to prove the following proposition.

\bp\label{prop_5}
Let $C_{k,n}$ denote the number of vertex-labeled normal networks with $n$ vertices and $k$ reticulation vertices which arise from the above caterpillar-skeleton. If $n$ is even then $C_{k,n}$ is zero, otherwise there is a positive constant $\tilde{d}_k$ such that
\[
C_{k,n}=n![z^n]C_k(z)\sim \tilde{d}_k\left(\frac{\sqrt{2}}{e}\right)^n n^{n+2k-1},
\]
as $n\rightarrow\infty$.
\ep

\bpf As in the proof of Corollary \ref{main-thm}, the asymptotic formula follows from (\ref{gen-func-cat}), where
\[
\tilde{d}_k=\frac{2\sqrt{2\pi}\tilde{e}_k(1/2)}{4^k\Gamma(2k-1/2)}.
\]

For the positivity claim, we will show that $\tilde{e}_k(1/2)$ is non-decreasing in $k$ from which
the claim follows by our result for $k=1$. In order to prove this, consider the
caterpillar-skeleton above with $k$ leaves. Denote the path consisting of the edges $e_1$ and
$e_2$ by $P$. Then, a subset of all normal networks generated by this caterpillar-skeleton of $k$
leaves is formed by normal networks which are generated by a caterpillar-skeleton with $k-1$
leaves to which a normal network with one reticulation vertex generated by $P$ is added. More
precisely, for the latter networks $g$ is connected to one of the subtrees attached to $e_1$ or
$e_2$ (such that the normal condition is satisfied), i.e., these networks arise from
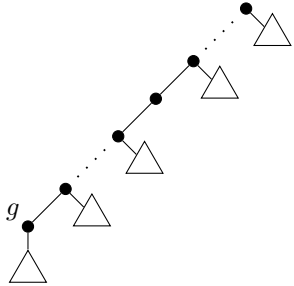
\begin{center}
\begin{tikzpicture}[scale=0.5,
		level/.style={thick},
		virtual/.style={thick,densely dashed},
		trans/.style={thick,<->,shorten >=2pt,shorten <=2pt,>=stealth},
		classical/.style={thin,double,<->,shorten >=4pt,shorten <=4pt,>=stealth}]

\draw (-11.2cm,-2.1cm) node[inner sep=1.5pt,circle,draw,fill] (1) {};
        \draw (-11.6cm,-1.7cm) node[inner sep=1.5pt,circle] (13) {$g$};
		\draw (-10.2cm,-1.1cm) node[inner sep=1.5pt,circle,draw,fill] (2) {};
		\draw (-8.8cm,0.3cm) node[inner sep=1.5pt,circle,draw,fill] (3) {};
		\draw (-7.8cm,1.3cm) node[inner sep=1.5pt,circle,draw,fill] (4) {};
		\draw (-6.8cm,2.3cm) node[inner sep=1.5pt,circle,draw,fill] (5) {};
		\draw (-5.4cm,3.7cm) node[inner sep=1.5pt,circle,draw,fill] (6) {};
        \draw (-11.2cm,-3.3cm) node[regular polygon,regular polygon sides=3,draw,inner sep=0.1cm] (7) {};
		\draw (-9.5cm,-1.8cm) node[regular polygon,regular polygon sides=3,draw,inner sep=0.1cm] (8) {};
		\draw (-8.1cm,-0.4cm) node[regular polygon,regular polygon sides=3,draw,inner sep=0.1cm] (9) {};
		\draw (-6.1cm,1.6cm) node[regular polygon,regular polygon sides=3,draw,inner sep=0.1cm] (11) {};
        \draw (-4.7cm,3cm) node[regular polygon,regular polygon sides=3,draw,inner sep=0.1cm] (12) {};
        \draw (-10cm,-0.9cm) node[inner sep=1.5pt,circle] (13) {};
        \draw (-9cm,0.1cm) node[inner sep=1.5pt,circle] (14) {};
        \draw (-6.6cm,2.5cm) node[inner sep=1.5pt,circle] (15) {};
        \draw (-5.6cm,3.5cm) node[inner sep=1.5pt,circle] (16) {};
		
		\draw (1)--(2);\draw (3)--(4);\draw (4)--(5);
        \draw (1)--(7);\draw (2)--(8);\draw (3)--(9);
        \draw (5)--(11);\draw (6)--(12);
		\draw[loosely dotted,line width=0.7pt] (13)--(14);
		\draw[loosely dotted,line width=0.7pt] (15)--(16);
\end{tikzpicture}
\end{center}
and are counted by
\begin{align*}
P(z)&=\partial_{y}\frac{z^{2}M(z,0)}{(1-z\tilde{M}(z,y))^2}\Big\vert_{y=0}
=\frac{8z^2-12z^4-(8z^2-4z^4)\sqrt{1-2z^2}}{(1-2z^2)^2},
\end{align*}
where $\tilde{M}(z,y)$ is as above. Consequently, the normal networks from the above mentioned subset are counted by
\[
C_{k-1}(z)P(z)=z\frac{\tilde{p}_k(z^2)-\tilde{q}_k(z^2)\sqrt{1-2z^2}}{(1-2z^2)^{2k-1/2}},
\]
where
\begin{align*}
\tilde{p}_k(z^2)&=(8z^2-12z^4)\tilde{e}_{k-1}(z^2)+(8z^2-4z^4)(1-2z^2)\tilde{f}_{k-1}(z^2);\\
\tilde{q}_k(z^2)&=(8z^2-4z^4)\tilde{e}_{k-1}(z^2)+(8z^2-12z^4)\tilde{f}_{k-1}(z^2).
\end{align*}
This gives, for odd $n$,
\[
n![z^{n}]C_{k-1}(z)P(z)\sim \tilde{g}_k\left(\frac{\sqrt{2}}{e}\right)^nn^{n+2k-1}
\]
with
\[
\tilde{g}_k=\frac{2\sqrt{2\pi}\tilde{e}_{k-1}(1/2)}{4^k\Gamma(2k-1/2)}>0.
\]
Moreover, since this counts a subclass of normal networks generated by a caterpillar-skeleton with
$k$ leaves, we have $\tilde{d}_k\geq \tilde{g}_k$ which gives $\tilde{e}_{k}(1/2)\geq \tilde{e}_{k-1}(1/2)$. This
proves our claim and thus the proposition is also proved.\epf

Finally, we would like to remark that in order to compute the multiplicative constant in the
asymptotic expression given in Corollary~\ref{main-thm} one has to understand the precise
structure of the generating functions for each Motzkin skeleton. Our investigations show that
the main contribution comes from the Motzkin skeletons for which the sparsened skeleton is a
(rooted, nonplane) tree with $k$ vertices. Since there is no explicit formula for the number of
such trees (but in fact there is an asymptotic solution; see \cite{AnaCombi}), we cannot expect to get
some explicit form for the constant, but only some expression
in terms of the number of rooted trees of size $k$. This observation may also be exploited to
derive upper bounds for the constant. With the help of Proposition~\ref{prop_5} lower bounds may
be derived as well. However, this needs some further investigations to understand the shape of the
polynomials $\tilde e_k(z)$ appearing in \eqref{gen-func-cat}.

\section{Counting Vertex-Labeled Tree-Child Networks}\label{vl-tc-net}

In this section, we will count (vertex-labeled) tree-child networks. As in the last section, we will first work out in detail the cases $k=1,2,3$, where, as for normal networks, we will show more precise results than stated in Theorem~\ref{main-thm-2}. The general case (and thus the proof of Theorem \ref{main-thm-2}) is then done in the last subsection below.

\subsection{Tree-child networks with one reticulation vertex}

We start with tree-child networks with one reticulation vertex which are again counted by using the Motzkin skeletons in Figure~\ref{fig_skeleton}.

\bp
The exponential generating function for vertex-labeled tree-child networks with
one reticulation vertex is
\begin{equation} \label{fctT1}
T_1(z)
=\frac{z^3\(1-\sqrt{1-2z^2}\,\)}{(1-2z^2)^{3/2}}=z\frac{\tilde{a}_1(z^2)-\tilde{b}_1(z^2)\sqrt{1-2z^2}}{(1-2z^2)^{3/2}},
\end{equation}
where
\[
\tilde{a}_1(z)=\tilde{b}_1(z)=z.
\]
\ep
\begin{proof}
We have to add an edge from $g$ in the Motzkin skeletons in Figure~\ref{fig_skeleton} which points to a unary (or red) vertex. Note that this edge is not allowed to point on a vertex on the path from $g$ to the root (since the resulting network must be a DAG), but is allowed to point to any vertex on the subtrees attached to these vertices. Moreover, the edge can also point to any non-root vertex in the subtree attached to $g$ (pointing on the root of this subtree is not allowed because we do not allow double edges).

This gives
\begin{align*}
T_1(z)&=\frac{z}{2}\partial_y \dfrac{\tilde{M}(z,y)}{1-z M(z,y)}\Big\vert_{y=0}
=\frac{z}{2}\left(\dfrac{M_y(z,0)-z^2-zM_b(z,0)}{1-z M(z,0)}+\dfrac{z M_{y}(z,0)M(z,0)}{(1-z M(z,0))^2}\right)
\end{align*}
where $\tilde{M}(z,y)$ is given in \eqref{tree-without-root}. Similar to the normal network case,
the factor $1/2$ compensates for the fact that each network is counted exactly twice by the above procedure. 
Now, by using \eqref{exp-M-y-b}
we obtain \eqref{fctT1}.
\end{proof}
From this, we obtain the following consequence; see the appendix for numerical data.
\begin{cor}
Let $T_{1,n}$ denote the number of vertex-labeled tree-child networks with $n$ vertices
and one reticulation vertex. If $n$ is even then $T_{1,n}$ is zero, otherwise
\[
T_{1,n}=n![z^n]T_1(z)=\(\frac{\sqrt2}{e}\)^n n^{n+1}\(\frac{\sqrt{2}}{2}-\frac{\sqrt{\pi}}{2}\cdot\frac{1}{\sqrt{n}}+{\mathcal O}\left(\frac{1}{n}\right)\),
\]
as $\nti$.
\end{cor}

\brem
Note that the constant of the second order term in the asymptotic expansion above is $-\sqrt{\pi}/2$ whereas that of the asymptotic expansion of $N_{1,n}$ is $-3\sqrt{\pi}/2$. Thus, the difference between normal networks and tree-child networks becomes visible only in the second order term (and the number of normal networks is of course smaller than the number of tree-child networks). The behavior for $k=2$ and $k=3$ is similar; see below.
\erem

\subsubsection*{Relationship to unicyclic networks revisited.} Again there is a close relationship to unicyclic networks and the alternative approach from Section \ref{rel-uni} can be used: either the root is in a cycle, but in which case now each vertex except the root can be the reticulation vertex, or the root is not in a cycle. This gives
\[
T_1(z)=zM(z,0)T_1(z)+\frac{1}{2}\sum_{\ell\geq 2}\ell z^{\ell+1}M(z,0)^{\ell}.
\]
Solving gives
\[
T_1(z)=\frac{\sum_{\ell\geq 2}\ell z^{\ell+1}M(z,0)^{\ell}}{2(1-zM(z,0))}=\frac{z^3M(z,0)^2(2-zM(z,0))}{2(1-zM(z,0))^3}.
\]
which by using the expression \eqref{tree-without-root} for $M(z,0)$ simplifies to (\ref{fctT1}).

\subsection{Tree-child networks with two reticulation vertices}

As for normal networks, the counting is done by using two variables $y_1$ and $y_2$ and the two types of Motzkin skeletons depicted in Figure~\ref{fig_skeleton_2}.

For trees attached to paths the situation is different from normal networks.
We never encounter different pointing rules between roots and internal vertices, but very well
between vertices on the path and vertices within the trees. Thus the red vertices in the third and fourth term on the right-hand side of the specification
for $\mathcal Q$, see \eqref{specQ}, fall into different classes of red vertices. In the third term, $\{\circ\}\times\mathcal Q\times (\{\bullet\}\times\tilde{\mathcal M})$, the red vertex is the root of the attached (red) tree. It can be treated like the red vertices within the tree and therfore
we mark it with $\tilde y$. A consequence of this is that we do not need to distinguish between red and white trees any more. Indeed, the second term of the
specification corresponds to having a white tree attached, the third one to having a red tree attached (to the path, in both cases). Since the red vertices fall
into the same class and are both marked by $\tilde y$, we may replace these two terms by one term corresponding to attaching simply a Motzkin tree.
The red vertex in the last term of \eqref{specQ} is on the path itself, thus marked by $y$. The other subtree cannot be a red tree by the tree-child condition.

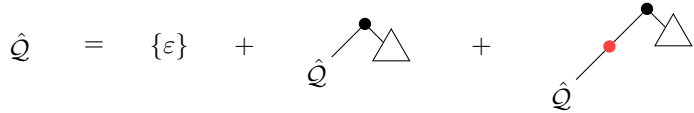
\begin{figure}[h]
\begin{center}
 \begin{tikzpicture}[
                scale=0.5,
                level/.style={thick},
                virtual/.style={thick,densely dashed},
                trans/.style={thick,<->,shorten >=2pt,shorten <=2pt,>=stealth},
                classical/.style={thin,double,<->,shorten >=4pt,shorten <=4pt,>=stealth}
                ]

                \draw (0cm,0.8cm) node[inner sep=1.5pt,circle] (2) {};
        \draw (-0.3cm,0.5cm) node[inner sep=1.5pt,circle] (1) {$\hat{\mathcal Q}$};
                \draw (1cm,1.8cm) node[inner sep=1.5pt,circle,draw,fill] (3) {};
                \draw (1.7cm,1.1cm) node[regular polygon,regular polygon sides=3,draw,inner sep=0.1cm] (4) {};

                \draw (2)--(3);\draw (3)--(4);

               \draw (6.5cm,0.2cm) node[inner sep=1.5pt,circle] (2) {};
        \draw (6.2cm,-0.1cm) node[inner sep=1.5pt,circle] (1) {$\hat{\mathcal Q}$};
                \draw (8.5cm,2.2cm) node[inner sep=1.5pt,circle,draw,fill] (3) {};
                \draw (9.2cm,1.5cm) node[regular polygon,regular polygon sides=3,draw,inner sep=0.1cm] (4) {};
        \draw (7.5cm,1.2cm) node[inner sep=1.5pt,circle,draw,fill,red!75] (5) {};

        \draw (2)--(5);\draw (5)--(3);\draw (3)--(4);

        \draw (4.1cm,1.2cm) node[inner sep=1.5pt,circle] (1) {$+$};
        \draw (-2.2cm,1.2cm) node[inner sep=1.5pt,circle] (1) {$+$};
        \draw (-4.2cm,1.2cm) node[inner sep=1.5pt,circle] (1) {$\{\varepsilon\}$};
        \draw (-6.2cm,1.2cm) node[inner sep=1.5pt,circle] (1) {$=$};
        \draw (-8.2cm,1.2cm) node[inner sep=1.5pt,circle] (1) {$\hat{\mathcal Q}$};
                \end{tikzpicture}
\end{center}	
\caption{The specification of the class $\hat{\mathcal Q}$ which is similar to that of ${\mathcal Q}$ ({\it cf}. Figure~\ref{QGF}) but with the second and third term merged. Also, now the subtree of the second term can be either red or white and that of the third term must be white. All the red vertices in these subtrees are counted by $\tilde{y}$; the other red vertices arising from the third term are counted by $y$.}	
\label{hat-Q}
\end{figure}

Altogether, this modification leads to a new class $\hat{\mathcal Q}$, specified by
\[
\hat{\mathcal Q}=\{\varepsilon\}\cup \{\circ\}\times \hat{\mathcal Q}\times \mathcal M \cup \{\circ\}
\times (\{\bullet\}\times\hat{\mathcal Q})\times \tilde{\mathcal M},
\]
see Figure~\ref{hat-Q}. We use this new class in \eqref{specP} instead of $\mathcal Q$ to specify the paths forming the basic building block for the Motzkin skeletons of tree-child networks. Call this new structure $\hat{\mathcal P}$. Then, we obtain the generating function
\[
\hat P(z,y,\tilde y,\hat y)=\frac{1+z\hat y}{1-z M(z,\tilde y)-z^2 y \tilde M(z,\tilde y)}.
\]
To summarize: The variable $y$ tells us which green vertex is allowed to point to vertices of the path (with the first vertex as possible exception),
$\tilde y$ which may point to vertices in the trees attached to the path, and $\hat y$ which may point to the first vertex of the path. We also make explicit a frequently appearing function:
\[
\hat P(z,0,\tilde y,0)=\rcp{1-zM(z,\tilde y)}.
\]

Now, the result for tree-child networks with two reticulation vertices is as follows.
\bp
The exponential generating function for vertex-labeled tree-child networks with
two reticulation vertices is
\begin{equation*}
T_2(z)=z\frac{\tilde{a}_2(z^2)-\tilde{b}_2(z^2)\sqrt{1-2z^2}}{(1-2z^2)^{7/2}},
\end{equation*}
where
\begin{equation*}
\tilde{a}_2(z)=-4z^4+\frac{21}{2}z^3-\frac{1}{2}z^2 \qquad\text{and}\qquad \tilde{b}_2(z)=9z^3-\frac{1}{2}z^2
\end{equation*}
\ep
\bpf
We start with the tree-child networks arising from the Motzkin skeletons on the left in Figure~\ref{fig_skeleton_2}. Here, $g_1$ and $g_2$ can point to all vertices in the attached subtrees except the root of the subtree attached to $g_1$. In addition, $g_2$ can also point to all vertices on the path between $g_1$ and $g_2$ except the vertex directly followed by $g_2$.

Overall, we obtain
\begin{align*}
T_{2,1}(z)&=
\partial_{y_1}\partial_{y_2}z^2\tilde M (z,y_1+y_2)\hat P(z,y_2,y_1+y_2,0)\hat P(z,0,y_1+y_2,0)\Big\vert_{y_1=0,y_2=0} \\
&=\partial_{y_1}\partial_{y_2}\frac{z^2\tilde M(z,y_1+y_2)}{(1-zM(z,y_1+y_2))(1-(z+z^2y_2)M(z,y_1+y_2)}\Big\vert_{y_1=0,y_2=0}.
\end{align*}

Now, consider the Motzkin skeletons on the right of Figure~\ref{fig_skeleton_2}. For the trees attached to the green vertices only pointing to the root is forbidden,
for all the other trees there is no pointing restriction. The analysis of the vertices on the paths is done path by path, as in the case of normal networks.
\begin{itemize}
\item Path $r$: No green vertex is allowed to point to the vertices of that path.
\item Path $k$: Poiting to all vertices is allowed for $g_2$, but $g_1$ may not point to that path at all. The situation for path $\ell$ is symmetric.
\end{itemize}

In this way, Motzkin skeletons which are not respecting the tree-child condition are generated as well: Indeed, $g_1$ may point to the
first vertex of $\ell$ and $g_2$ to the first vertex of $k$, such that the common ancestor $y$ has two red children. But for these networks, all paths
are of the same type, since there cannot be any red vertices in any subtree.
This gives
\begin{align*}
T_{2,2}(z)=&\rcp2\partial_{y_1}\partial_{y_2}\frac{z^3\tilde M(z,y_1+y_2)^2}{1-zM(z,y_1+y_2)}\hat P(z,y_2,y_1+y_2,y_2)
\hat P(z,y_1,y_1+y_2,y_1)\Big\vert_{y_1=0,y_2=0}
\\
&-\rcp2\frac{z^5M(z,0)^2}{(1-zM(z,0))^3}.
\end{align*}

The exponential generating function for vertex-labeled tree-child networks is now obtained as $T_2(z)=(T_{2,1}(z)+T_{2,2}(z))/4$.
Plugging in the above expressions and simplifying gives the result.\epf

As a consequence, we have the following; see the appendix for numerical data.

\begin{cor}
Let $T_{2,n}$ denote the number of vertex-labeled tree-child networks with $n$
vertices and two reticulation vertices. If $n$ is even then $T_{2,n}$ is zero, otherwise
\[
T_{2,n}=n![z^n]T_{2}(z)=\(\frac{\sqrt2}{e}\)^n n^{n+3}\(\frac{\sqrt{2}}{16}-\frac{\sqrt{\pi}}{8}\cdot\frac{1}{\sqrt{n}}+{\mathcal O}\left(\frac{1}{n}\right)\),
\]
as $\nti$.
\end{cor}

\subsection{Tree-child networks with three reticulation vertices}

In this case we use the four different types of Motzkin skeletons depicted in Figure~\ref{fig_skeleton_3}, Figure~\ref{fig_skeleton_4} and Figure~\ref{fig_skeleton_5}. Moreover,
we use
the $\Y$ operator from Section~\ref{normal-k=3}.

We start with the tree-child networks arising from the Motzkin skeletons depicted on the left of Figure~\ref{fig_skeleton_3}. The possibilities for the pointings of the edges starting at $g_1, g_2$ and $g_3$ are similar as in the first case for $k=2$ (see above). All these edges may target any non-root vertex in the tree attached to $g_1$ and any vertex in all the other trees. Concerning the vertices on the spine, we have some restrictions. The edge from $g_1$ may not end at any vertex from $\ell_1$, for the first vertex this applies even to $g_2$. Similarly, the edges from $g_1$ and $g_2$ may not point to any vertex of $\ell_2$, and
no green vertex may point to the first vertex of $\ell_2$ as well as to any vertex of $\ell_1$.

Overall, we obtain for this Motzkin skeleton
	\begin{align*}
		T_{3,1}(z)&=\Y\(\frac{z^3\tilde M(z,y_1+y_2+y_3)\hat P(z,y_3,y_1+y_2+y_3,0)\hat P(z,y_2+y_3,y_1+y_2+y_3,0)}{1-zM(z,y_1+y_2+y_3)}\).
		\end{align*}		

For the other cases, a similar reasoning for the possible pointings of the edges starting from $g_1, g_2$ and $g_3$ can be used. Furthermore, we have pay attention
to the Motzkin skeletons we generate which are not tree-child. These cases are those where two green vertices point to the children of a latest common ancestor of two
green vertices, but the third green vertex has some freedom in pointing. We refrain from giving details and just list the obtained expressions. The reader
is invited to derive them herself.

For the Motzkin skeletons on the right of Figure~\ref{fig_skeleton_3}, we obtain
		\begin{align*}
			T_{3,2}(z)=&\rcp2
			\Y\left(\frac{z^4\tilde M(z,y_1+y_2+y_3)^2}{1-zM(z,y_1+y_2+y_3)}\hat P(z,y_3,y_1+y_2+y_3,0)\right.
			\\
			&\qquad\times \hat P(z,y_1+y_3,y_1+y_2+y_3,y_1+y_3) \hat P(z,y_2+y_3,y_1+y_2+y_3,y_2+y_3)\Bigg)
			\\
			&-\rcp2\partial_{y_1}\frac{z^6\tilde M(z,y_1)^2}{(1-zM(z,y_1))^3}\hat P(z,y_1,y_1,0)\Big\vert_{y_1=0}
			-\rcp2\partial_{y_2}\frac{z^6\tilde M(z,y_2)^2}{(1-zM(z,y_2))^3}\hat P(z,y_2,y_2,0)\Big\vert_{y_2=0}
			\\
			&-\rcp2\partial_{y_3}\frac{z^6\tilde M(z,y_3)^2}{1-zM(z,y_3)}\hat P(z,y_3,y_3,0)^3\Big\vert_{y_3=0}.
		\end{align*}

For the Motzkin skeletons depicted in Figure~\ref{fig_skeleton_4}, we obtain
\begin{align*}
		T_{3,3}(z)=&
		\Y\left(\frac{z^4\tilde M(z,y_1+y_2+y_3)^2}{1-zM(z,y_1+y_2+y_3)}\hat P(z,y_2+y_3,y_1+y_2+y_3,y_2+y_3)\right.
		\\
		&\ \quad\times \hat P(z,y_1,y_1+y_2+y_3,y_1) \hat P(z,y_1+y_2,y_1+y_2+y_3,0)\Bigg) \\
		&-\partial_{y_2}\frac{z^6\tilde M(z,y_2)^2}{(1-zM(z,y_2))^2}\hat P(z,y_2,y_2,0)^2\Big\vert_{y_2=0}
		-\partial_{y_3}\frac{z^6\tilde M(z,y_3)^2}{(1-zM(z,y_3))^3}\hat P(z,y_3,y_3,0)\Big\vert_{y_3=0}.
		\end{align*}

For the final case, consider the Motzkin skeletons depicted in Figure~\ref{fig_skeleton_5}. Here, the generating function is given by
\begin{align*}
T_{3,4}(z)=&\rcp2
\Y\left(\frac{z^5\tilde M(z,y_1+y_2+y_3)^3}{1-zM(z,y_1+y_2+y_3)}
\right.
\\
&\qquad\times \hat P(z,y_1+y_2,y_1+y_2+y_3,y_1+y_2)\hat P(z,y_1+y_3,y_1+y_2+y_3,y_1+y_3)
\\
&\qquad\times \hat P(z,y_2+y_3,y_1+y_2+y_3,y_2+y_3)\hat P(z,y_3,y_1+y_2+y_3,y_3)\Bigg)\\
&-\rcp2\partial_{y_3}\frac{z^7\tilde M(z,y_3)^3}{(1-zM(z,y_3))^2}\hat P(z,y_3,y_3,y_3)\hat P(z,y_3,y_3,0)^2\Big\vert_{y_3=0}
\\
&-\partial_{y_2}\frac{z^7\tilde M(z,y_2)^3}{(1-zM(z,y_2))^3}\hat P(z,y_2,y_2,y_2)\hat P(z,y_2,y_2,0)\Big\vert_{y_2=0}
\\
&-\partial_{y_1}\frac{z^7\tilde M(z,y_1)^3}{(1-zM(z,y_1))^3}\hat P(z,y_1,y_1,y_1)\hat P(z,y_1,y_1,0)\Big\vert_{y_1=0}.
\end{align*}

The exponential generating function for vertex-labeled tree-child networks is obtained as $T_{3}(z)=(T_{3,1}(z)+T_{3,2}(z)+T_{3,3}(z)+T_{3,4}(z))/8$
after all. This gives the following result.				

\bp
The exponential generating function for vertex-labeled tree-child networks with
three reticulation vertices is
\begin{equation*}
T_3(z)=z\frac{\tilde{a}_3(z^2)-\tilde{b}_3(z^2)\sqrt{1-2z^2}}{(1-2z^2)^{11/2}},
\end{equation*}
where
\[
\tilde{a}_3(z)=-6z^6+\frac{249}{2}z^5+20 z^4-z^3
\]
and
\[
\tilde{b}_3(z)=30z^6+144z^5+19z^4-z^3
\]
\ep

As a consequence, we have the following result; see the appendix for numerical data.

\begin{cor}
Let $T_{3,n}$ denote the number of vertex-labeled tree-child networks with $n$ vertices
and three reticulation vertices. If $n$ is even then $T_{3,n}$ is zero, otherwise
\[
T_{3,n}=n![z^n]T_{3}(z)=\(\frac{\sqrt2}{e}\)^n n^{n+5}\(\frac{\sqrt{2}}{192}-\frac{\sqrt{\pi}}{64}\cdot\frac{1}{\sqrt{n}}+{\mathcal O}\left(\frac{1}{n}\right)\),
\]
as $\nti$.
\end{cor}

\subsection{Tree-child networks with a fixed number of reticulation vertices}

In this subsection, we will prove Theorem \ref{main-thm-2} which is deduced from the following proposition.
\bp
For the numbers of vertex-labeled normal networks $N_{k,n}$ and vertex-labeled tree-child networks $T_{k,n}$,
\[
T_{k,n}=N_{k,n}\left(1+{\mathcal O}\left(\frac{1}{\sqrt{n}}\right)\right), \text{ as } n\to\infty.
\]
\ep

\bpf First, observe that $T_{k,n}-N_{k,n}$ is bounded by the number of networks which arise from all types of Motzkin skeletons where for each green vertex we consider all possibilities of adding an edge such that the normal condition is violated (note that this is an over-estimate of the difference). Thus, we only have to count the number of such networks which arise from a fixed type of Motzkin skeletons and a fixed green vertex. Similar to the proof of Proposition \ref{gen-normal-k}, the largest number will come from the Motzkin skeletons where the green vertices are the leaves (this will become clear by applying the same arguments as below to all other Motzkin skeletons).

Now, fix such a type of Motzkin skeletons and one of its green vertices. Then, for this vertex, we will have the following options.
\begin{itemize}
\item The green vertex points to one of the subtrees attached to the leaves of the skeletons. For the exponential generating function this gives
\begin{equation*}
\partial_{y_2}\cdots\partial_{y_k}\frac{G'_1(z,y)\cdots G_{s}(z,y)}{(1-G_{s+1}(z,y))\cdots(1-G_{s+2k-1}(z,y))}\Big\vert_{y_2=0,\ldots,y_k=0},
\end{equation*}
where the derivative comes from choosing a vertex in the subtree as end point of the green vertex. (Here, and below $y$ is the sum of $y_i$'s with $2\leq i\leq k$ and not all of the $y_i$'s must be present; also which are present can differ from one occurrence to the next.)
\item The green vertex points to the root of a Motzkin tree from ${\mathcal M}$ attached to the path from the green vertex to the root or attached to some of the edges of the sparsened skeleton on a path from the green vertex to a leaf. Then, we have
\begin{equation*}
\partial_{y_2}\cdots\partial_{y_k}\frac{G_1(z,y)\cdots G_{s}(z,y)}{(1-G_{s+1}(z,y))\cdots(1-G_{s+2k-1}(z,y))(1-G_{s+2k}(z,y))}\Big\vert_{y_2=0,\ldots,y_k=0},
\end{equation*}
where the additional term comes from the fact that now one edge was split into two edges by the above pointing.
\item The green vertex points to the first vertex on one of the branches attached to the path from the green vertex to the root. Then, we have
\begin{equation*}
\partial_{y_2}\cdots\partial_{y_k}\frac{G_1(z,y)\cdots G_{s}(z,y)}{(1-G_{s+1}(z,y))\cdots(1-G_{s+2k-1}(z,y))}\Big\vert_{y_2=0,\ldots,y_k=0}.
\end{equation*}
\end{itemize}

The exponential generating function of all networks arising from these Motzkin skeletons and the green vertex are a sum of generating functions of the above three types. Thus, from Lemma \ref{tech-lmm}, we obtain that this generating function has the form
\[
\frac{c(z)-d(z)\sqrt{1-2z^2}}{(1-2z^2)^{p}},
\]
where $c(z)$ and $d(z)$ are suitable polynomials and the maximum of $p$ is as follows: note that without the derivatives in the above expressions, $p$ would be at most $k$ (this bound is taken on in the first two cases, but not in the last case where $p$ is at most $k-1/2$); also, because of Lemma \ref{tech-lmm}, each derivative increases this bound by one. Thus, $p$ is at most $2k-1$.

Now, with the same arguments as in the proof of Corollary \ref{cor-gen}, we obtain that the exponential generating function of the above number has the form
\[
z\frac{\tilde{c}(z^2)-\tilde{d}(z^2)\sqrt{1-2z^2}}{(1-2z^2)^{2k-1}},
\]
where $\tilde{c}(z)$ and $\tilde{d}(z)$ are suitable polynomials. Singularity analysis gives then the bound
\[
{\mathcal O}\left(\left(\frac{\sqrt{2}}{e}\right)^nn^{n+2k-3/2}\right).
\]
Summing over all possible type of Motzkin skeletons and all green vertices, we obtain the same bound for $T_{k,n}-N_{k,n}$ which proves the claimed result.\epf

\section{Counting Leaf-Labeled Normal and Tree-Child Networks}\label{ll-net}

In this section, we will count leaf-labeled normal and tree-child networks with $\ell$ leaves and $k$ reticulation vertices (recall that we denoted their numbers by $\tilde{N}_{k,\ell}$ and $\tilde{T}_{k,\ell}$, respectively). The counting results will follow from those for vertex-labeled networks since there is a close relationship between leaf-labeled normal and tree-child networks and vertex-labeled ones. To see this, we need to recall two lemmas from \cite{MSW15}.

\begin{lemma}[see \cite{MSW15}]
For any phylogenetic network with $\ell$ leaves, $k$ reticulation vertices and $n$ vertices, we have
\[
\ell+k=\frac{n+1}{2}.
\]
(Recall that $n$ is always odd.)
\end{lemma}
\begin{lemma}[see \cite{MSW15}]\label{des-equal}
The descendant sets for any two non-leaf vertices in a tree-child network (and thus also normal network) are different.
\end{lemma}

These two lemmas immediately imply that
\[
N_{k,2\ell+2k-1}=\binom{2\ell+2k-1}{\ell}(\ell+2k-1)!\tilde{N}_{k,\ell}.
\]

To see this, note that all vertex-labeled normal networks with $2\ell+2k-1$ vertices and $k$ reticulation vertices can be constructed as follows: start with a (fixed) leaf-labeled normal network with $\ell$ leaves and $k$ reticulation vertices. Then, choose $\ell$ labels from the set $2\ell+2k-1$ labels and re-label the leaves of the fixed network such that the order is preserved. Finally, label the remaining $\ell+2k-1$ vertices by any permutation of the set of remaining $\ell+2k-1$ labels. By the above two lemmas, in this way every vertex-labeled normal network is obtained exactly once.

The above now implies that
\[
\tilde{N}_{k,\ell}=\frac{\ell!}{(2\ell+2k-1)!}N_{k,2\ell+2k-1}
\]
from which an asymptotic result follows by Theorem~\ref{main-thm-1} and Stirling's formula. Similarly, an asymptotic result for leaf-labeled tree-child networks is obtained from
\[
\tilde{T}_{k,\ell}=\frac{\ell!}{(2\ell+2k-1)!}T_{k,2\ell+2k-1}.
\]
Overall, we obtain the following theorem.

\begin{theo}
For the numbers $\tilde{N}_{k,\ell}$ and $\tilde{T}_{k,\ell}$ of leaf-labeled normal networks resp. leaf-labeled tree-child networks with $k\geq 1$ reticulation vertices, we have
\[
\tilde{N}_{k,\ell}\sim\tilde{T}_{k,\ell}\sim2^{3k-1}c_k\left(\frac{2}{e}\right)^{\ell}\ell^{\ell+2k-1},\qquad(\ell\rightarrow\infty)
\]
where $c_k$ is as in Theorem \ref{main-thm-1}.
\end{theo}

\section{Conclusion}\label{con}
In this paper, we considered the counting problem of phylogenetic networks which is largely
unsolved. We devised an approach, based on generating functions and analytic combinatorics, to
solve this problem for two important subclasses of phylogenetic networks, namely, tree-child and
normal networks, provided that the number of reticulation vertices is fixed as the size of the
network tends to infinity. The latter restriction is necessary for our method to work. Indeed, the
combinatorial setup we developed in this paper is the construction of a sequence of combinatorial
classes (for each given number of reticulation vertices, we contruct a separate class). The actual
distribution of the reticulation vertices is then -- on the level of generating functions --
realized by differentiations. Letting $k$ tend to infinity, when $n$ tends to infinity, means that
we have to cope with a growing number of differentiations and it is not clear how this changes the
qualitative nature of the generating function. We certainly cannot expect that $N_k(z)$ keeps the
shape \eqref{gen-form} when $k$ depends on $n$ and gets large with growing $n$. Thus, we have
to leave the question of counting phylogenetic networks when $k$ is allowed to grow with $n$ open.

Apart from this, the most obvious other question about the results of this paper is the following: why is our method only applied to
subclasses of phylogenetic networks? In fact, our method can probably be extended to count general
networks with a fixed number of reticulation vertices, too, but for this some further work has to
be done. We will explain now why.

First, for vertex-labeled networks our method above relied on the use of Motzkin skeletons, which have green and red vertices, and all of them are unary vertices. Recall that these vertices arise by deleting an edge for each reticulation vertex which was colored red (the green vertices are then the other endpoints of the deleted edges). If one considers general phylogenetic networks, then the colored vertices in the Motzkin skeleton can be leaves as well. In order to see this consider the following networks (which show all possible types of how such leaves can occur):

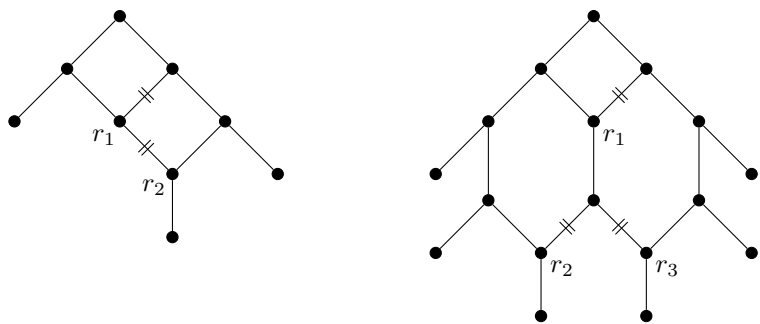
\begin{center}
\begin{tikzpicture}[scale=0.7]
\draw (0cm,0cm) node[inner sep=1.5pt,circle,draw,fill] (1) {};
\draw (-1cm,-1cm) node[inner sep=1.5pt,circle,draw,fill] (2) {};
\draw (1cm,-1cm) node[inner sep=1.5pt,circle,draw,fill] (3) {};
\draw (-2cm,-2cm) node[inner sep=1.5pt,circle,draw,fill] (4) {};
\draw (0cm,-2cm) node[inner sep=1.5pt,circle,draw,fill] (5) {};
\draw (2cm,-2cm) node[inner sep=1.5pt,circle,draw,fill] (6) {};
\draw (1cm,-3cm) node[inner sep=1.5pt,circle,draw,fill] (7) {};
\draw (3cm,-3cm) node[inner sep=1.5pt,circle,draw,fill] (8) {};
\draw (1cm,-4.2cm) node[inner sep=1.5pt,circle,draw,fill] (9) {};
\draw (-0.3cm,-2.3cm) node[inner sep=1.5pt,circle] (10){$r_1$};
\draw (0.65cm,-3.25cm) node[inner sep=1.5pt,circle] (11) {$r_2$};
\draw (0.27cm,-1.33cm) node[inner sep=1.5pt,circle] (12) {};
\draw (0.67cm,-1.73cm) node[inner sep=1.5pt,circle] (13) {};
\draw (0.33cm,-1.27cm) node[inner sep=1.5pt,circle] (14) {};
\draw (0.73cm,-1.67cm) node[inner sep=1.5pt,circle] (15) {};
\draw (0.67cm,-2.27cm) node[inner sep=1.5pt,circle] (16) {};
\draw (0.27cm,-2.67cm) node[inner sep=1.5pt,circle] (17) {};
\draw (0.73cm,-2.33cm) node[inner sep=1.5pt,circle] (18) {};
\draw (0.33cm,-2.73cm) node[inner sep=1.5pt,circle] (19) {};

\draw (1)--(2);\draw (1)--(3);\draw (2)--(4);\draw (2)--(5);
\draw (3)--(5);\draw (3)--(6);\draw (5)--(7);\draw (6)--(8);
\draw (7)--(9);\draw (6)--(7);\draw (12)--(13);\draw (14)--(15);
\draw (16)--(17);\draw (18)--(19);

\draw (9cm,0cm) node[inner sep=1.5pt,circle,draw,fill] (1) {};
\draw (8cm,-1cm) node[inner sep=1.5pt,circle,draw,fill] (2) {};
\draw (10cm,-1cm) node[inner sep=1.5pt,circle,draw,fill] (3) {};
\draw (7cm,-2cm) node[inner sep=1.5pt,circle,draw,fill] (4) {};
\draw (9cm,-2cm) node[inner sep=1.5pt,circle,draw,fill] (5) {};
\draw (11cm,-2cm) node[inner sep=1.5pt,circle,draw,fill] (6) {};
\draw (6cm,-3cm) node[inner sep=1.5pt,circle,draw,fill] (7) {};
\draw (12cm,-3cm) node[inner sep=1.5pt,circle,draw,fill] (8) {};
\draw (7cm,-3.5cm) node[inner sep=1.5pt,circle,draw,fill] (9) {};
\draw (9cm,-3.5cm) node[inner sep=1.5pt,circle,draw,fill] (10) {};
\draw (11cm,-3.5cm) node[inner sep=1.5pt,circle,draw,fill] (11) {};
\draw (6cm,-4.5cm) node[inner sep=1.5pt,circle,draw,fill] (12) {};
\draw (8cm,-4.5cm) node[inner sep=1.5pt,circle,draw,fill] (13) {};
\draw (10cm,-4.5cm) node[inner sep=1.5pt,circle,draw,fill] (14) {};
\draw (12cm,-4.5cm) node[inner sep=1.5pt,circle,draw,fill] (15) {};
\draw (8cm,-5.7cm) node[inner sep=1.5pt,circle,draw,fill] (16) {};
\draw (10cm,-5.7cm) node[inner sep=1.5pt,circle,draw,fill] (17) {};
\draw (9.4cm,-2.3cm) node[inner sep=1.5pt,circle] (18) {$r_1$};
\draw (8.4cm,-4.8cm) node[inner sep=1.5pt,circle] (19) {$r_2$};
\draw (10.4cm,-4.8cm) node[inner sep=1.5pt,circle] (20) {$r_3$};
\draw (9.27cm,-1.33cm) node[inner sep=1.5pt,circle] (21) {};
\draw (9.67cm,-1.73cm) node[inner sep=1.5pt,circle] (22) {};
\draw (9.33cm,-1.27cm) node[inner sep=1.5pt,circle] (23) {};
\draw (9.73cm,-1.67cm) node[inner sep=1.5pt,circle] (24) {};
\draw (8.27cm,-3.83cm) node[inner sep=1.5pt,circle] (25) {};
\draw (8.67cm,-4.23cm) node[inner sep=1.5pt,circle] (26) {};
\draw (8.33cm,-3.77cm) node[inner sep=1.5pt,circle] (27) {};
\draw (8.73cm,-4.17cm) node[inner sep=1.5pt,circle] (28) {};
\draw (9.67cm,-3.77cm) node[inner sep=1.5pt,circle] (29) {};
\draw (9.27cm,-4.17cm) node[inner sep=1.5pt,circle] (30) {};
\draw (9.73cm,-3.83cm) node[inner sep=1.5pt,circle] (31) {};
\draw (9.33cm,-4.23cm) node[inner sep=1.5pt,circle] (32) {};

\draw (1)--(2);\draw (1)--(3);\draw (2)--(4);\draw (2)--(5);
\draw (3)--(5);\draw (3)--(6);\draw (4)--(7);\draw (6)--(8);
\draw (4)--(9);\draw (5)--(10);\draw (6)--(11);\draw (9)--(12);
\draw (9)--(13);\draw (10)--(13);\draw (10)--(14);\draw (11)--(14);
\draw (11)--(15);\draw (13)--(16);\draw (14)--(17);\draw (21)--(22);
\draw (23)--(24);\draw (25)--(26);\draw (27)--(28);\draw (29)--(30);
\draw (31)--(32);
\end{tikzpicture}
\end{center}

In the network on the left, if the indicated edges are deleted, then $r_1$ becomes a leaf (which is colored both green and red). On the other hand, in the network on the right, after deleting the indicated edges, the vertex which was connected to $r_2$ and $r_3$ becomes a leaf (which is colored green). Also, note that in the second cases the number of green and red vertices in the resulting Motzkin skeleton is not the same (unless one considers the leaf in the second case to be colored ``double-green"). So, in order to consider the counting problem for vertex-labeled general networks, more possibilities for the Motzkin skeletons must be considered.

However, we in fact suspect that all the above mentioned additional possibilities for the Motzkin
skeletons are asymptotically negligible since they lead to restrictions: In the first case one
green vertex must be connected to the red-green leaf which reduces the number of differentiations
in the expression for the exponential generating function by one, and in the second case the
number of green vertices and thus the number of edges in the Motzkin skeleton is reduced by one
(which also leads to a contribution of smaller order; e.g. see the proof of
Proposition~\ref{gen-normal-k}). Moreover, for the general networks arising from the Motzkin
skeletons defined in this paper, we also expect that those networks which do not satisfy the
tree-child property are rare, because again, when the tree-child property is not satisfied, then
one has severe restrictions. So, overall, we guess that also for vertex-labeled general
phylogenetic networks the same first-order asymptotics as in Theorem~\ref{main-thm-1} and
Theorem~\ref{main-thm-2} holds.

Second, the method for the counting of leaf-labeled general phylogenetic networks with $k$
reticulation vertices will also be different from the one used in this paper since we do not have
a simple connection between the vertex-labeled and the leaf-labeled case anymore. (This is because
Lemma~\ref{des-equal} no longer holds.) Thus, one has to cope with symmetries. However, it is
expected that phylogenetic networks which have vertices with equal sets of descendants are rare
and thus one again expects the same first-order asymptotics as in Theorem~\ref{main-thm}. We may
come back to these questions elsewhere.

Finally, we mention that tree-child and normal phylogenetic networks are special classes of
directed acyclic graphs. Though their enumeration is in general not easy, there is already a
fairly rich basis of enumeration results on general directed acyclic graphs \cite{BRRW86, BeRo88,
Ro76, Ro71, Ro92} as well as more sophisticated studies on their shape \cite{Li75,McK89} which
can probably be exploited in order to extend our results or get a finer analysis on the structure
of random phylogenetic networks.

\section*{Acknowledgment} We thank the reviewers for many helpful comments.

\bibliographystyle{plain}
\bibliography{phylo-5}

\begin{thebibliography}{10}

\bibitem{AlAl16}
Nikita Alexeev and Max~A. Alekseyev.
\newblock Combinatorial scoring of phylogenetic networks.
\newblock In {\em Computing and combinatorics}, volume 9797 of {\em Lecture
  Notes in Comput. Sci.}, pages 560--572. Springer, [Cham], 2016.

\bibitem{BRRW86}
Edward~A. Bender, L.~Bruce Richmond, Robert~W. Robinson, and Nick~C. Wormald.
\newblock The asymptotic number of acyclic digraphs. {I}.
\newblock {\em Combinatorica}, 6(1):15--22, 1986.

\bibitem{BeRo88}
Edward~A. Bender and Robert~W. Robinson.
\newblock The asymptotic number of acyclic digraphs. {II}.
\newblock {\em J. Combin. Theory Ser. B}, 44(3):363--369, 1988.

\bibitem{Bo16}
Mikl\'os B\'ona.
\newblock On the number of vertices of each rank in {$k$}-phylogenetic trees.
\newblock {\em Discrete Math. Theor. Comput. Sci.}, 18(3):Paper No. 7, 7, 2016.

\bibitem{BoFl09}
Mikl\'os B\'ona and Philippe Flajolet.
\newblock Isomorphism and symmetries in random phylogenetic trees.
\newblock {\em J. Appl. Probab.}, 46(4):1005--1019, 2009.

\bibitem{BoLiSe17}
Magnus Bordewich, Simone Linz, and Charles Semple.
\newblock Lost in space? {G}eneralising subtree prune and regraft to spaces of
  phylogenetic networks.
\newblock {\em J. Theoret. Biol.}, 423:1--12, 2017.

\bibitem{BoSe16a}
Magnus Bordewich and Charles Semple.
\newblock Determining phylogenetic networks from inter-taxa distances.
\newblock {\em J. Math. Biol.}, 73(2):283--303, 2016.

\bibitem{BoSe16b}
Magnus Bordewich and Charles Semple.
\newblock Reticulation-visible networks.
\newblock {\em Adv. in Appl. Math.}, 78:114--141, 2016.

\bibitem{CaRoVa}
Gabriel Cardona, Francesc Rossello, and Gabriel Valiente.
\newblock Comparison of tree-child phylogenetic networks.
\newblock {\em IEEE/ACM Trans. Comput. Biol. Bioinform.}, 6:552--569, 2009.

\bibitem{ChWa12}
Zhi-Zhong Chen and Lusheng Wang.
\newblock Algorithms for reticulate networks of multiple phylogenetic trees.
\newblock {\em IEEE/ACM Transactions on Computational Biology and
  Bioinformatics}, 9(2):372--384, 2012.

\bibitem{CLS14}
Paul Cordue, Simone Linz, and Charles Semple.
\newblock Phylogenetic networks that display a tree twice.
\newblock {\em Bull. Math. Biol.}, 76(10):2664--2679, 2014.

\bibitem{CEJM13}
\'{E}va Czabarka, P\'{e}ter~L. Erd\H{o}s, Virginia Johnson, and Vincent
  Moulton.
\newblock Generating functions for multi-labeled trees.
\newblock {\em Discrete Appl. Math.}, 161(1-2):107--117, 2013.

\bibitem{DiRo17}
Filippo Disanto and Noah~A. Rosenberg.
\newblock Enumeration of ancestral configurations for matching gene trees and
  species trees.
\newblock {\em J. Comput. Biol.}, 24(9):831--850, 2017.

\bibitem{FlOd82}
Philippe Flajolet and Andrew Odlyzko.
\newblock The average height of binary trees and other simple trees.
\newblock {\em J. Comput. System Sci.}, 25(2):171--213, 1982.

\bibitem{FO90}
Philippe Flajolet and Andrew Odlyzko.
\newblock Singularity analysis of generating functions.
\newblock {\em SIAM J. Discrete Math.}, 3(2):216--240, 1990.

\bibitem{FlPr86}
Philippe Flajolet and Helmut Prodinger.
\newblock Register allocation for unary-binary trees.
\newblock {\em SIAM J. Comput.}, 15(3):629--640, 1986.

\bibitem{AnaCombi}
Philippe Flajolet and Robert Sedgewick.
\newblock {\em Analytic Combinatorics}.
\newblock Cambridge University Press, Cambridge, 2009.

\bibitem{FlSt80}
Philippe Flajolet and Jean-Marc Steyaert.
\newblock On the analysis of tree-matching algorithms.
\newblock In {\em Automata, languages and programming ({P}roc. {S}eventh
  {I}nternat. {C}olloq., {N}oordwijkerhout, 1980)}, volume~85 of {\em Lecture
  Notes in Comput. Sci.}, pages 208--219. Springer, Berlin-New York, 1980.

\bibitem{FoRo80}
Leslie~R. Foulds and Robert~W. Robinson.
\newblock Determining the asymptotic number of phylogenetic trees.
\newblock In {\em Combinatorial mathematics, {VII} ({P}roc. {S}eventh
  {A}ustralian {C}onf., {U}niv. {N}ewcastle, {N}ewcastle, 1979)}, volume 829 of
  {\em Lecture Notes in Math.}, pages 110--126. Springer, Berlin, 1980.

\bibitem{FoRo88}
Leslie~R. Foulds and Robert~W. Robinson.
\newblock Enumerating phylogenetic trees with multiple labels.
\newblock In {\em Proceedings of the {F}irst {J}apan {C}onference on {G}raph
  {T}heory and {A}pplications ({H}akone, 1986)}, volume~72, pages 129--139,
  1988.

\bibitem{FrSt15}
Andrew~R. Francis and Mike Steel.
\newblock Tree-like reticulation networks---when do tree-like distances also
  support reticulate evolution?
\newblock {\em Math. Biosci.}, 259:12--19, 2015.

\bibitem{KeSc14}
Steven Kelk and Celine Scornavacca.
\newblock Constructing minimal phylogenetic networks from softwired clusters is
  fixed parameter tractable.
\newblock {\em Algorithmica}, 68(4):886--915, 2014.

\bibitem{LSS13}
Simone Linz, Katherine St.~John, and Charles Semple.
\newblock Counting trees in a phylogenetic network is \#{P}-complete.
\newblock {\em SIAM J. Comput.}, 42(4):1768--1776, 2013.

\bibitem{Li75}
Valery~A. Liskovec.
\newblock The number of maximal vertices of a random acyclic digraph.
\newblock {\em Teor. Verojatnost. i Primenen.}, 20(2):412--421, 1975.

\bibitem{MSW15}
Colin McDiarmid, Charles Semple, and Dominic Welsh.
\newblock Counting phylogenetic networks.
\newblock {\em Ann. Comb.}, 19(1):205--224, 2015.

\bibitem{McK89}
Brendan~D. McKay.
\newblock On the shape of a random acyclic digraph.
\newblock {\em Math. Proc. Cambridge Philos. Soc.}, 106(3):459--465, 1989.

\bibitem{Ro71}
Robert~W. Robinson.
\newblock Counting labeled acyclic digraphs.
\newblock pages 239--273, 1973.

\bibitem{Ro76}
Robert~W. Robinson.
\newblock Counting unlabeled acyclic digraphs.
\newblock pages 28--43. Lecture Notes in Math., Vol. 622, 1977.

\bibitem{Ro92}
Vitalii~I. Rodionov.
\newblock On the number of labeled acyclic digraphs.
\newblock {\em Discrete Math.}, 105(1-3):319--321, 1992.

\bibitem{Ro07}
Noah~A. Rosenberg.
\newblock Counting coalescent histories.
\newblock {\em J. Comput. Biol.}, 14(3):360--377, 2007.

\bibitem{Sc}
Ernst Schr\"oder.
\newblock Vier kombinatorische {P}robleme.
\newblock {\em Z. Math. Phys.}, 15:361--376, 1870.

\bibitem{Se16}
Charles Semple.
\newblock Phylogenetic networks with every embedded phylogenetic tree a base
  tree.
\newblock {\em Bull. Math. Biol.}, 78(1):132--137, 2016.

\bibitem{Se17}
Charles Semple.
\newblock Size of a phylogenetic network.
\newblock {\em Discrete Appl. Math.}, 217(part 2):362--367, 2017.

\bibitem{SeSt06}
Charles Semple and Mike Steel.
\newblock Unicyclic networks: compatibility and enumeration.
\newblock {\em IEEE/ACM Trans. Comput. Biology Bioinform.}, 3:84--91, 2006.

\bibitem{vKS16}
Leo van Iersel, Steven Kelk, and Celine Scornavacca.
\newblock Kernelizations for the hybridization number problem on multiple
  nonbinary trees.
\newblock {\em J. Comput. System Sci.}, 82(6):1075--1089, 2016.

\bibitem{Wi}
Stephen~J. Willson.
\newblock Properties of normal phylogenetic networks.
\newblock {\em Bull. Math. Biol.}, 72(2):340--358, 2010.

\end{thebibliography}

\newpage
\appendix
\section{Number of Networks and Their Asymptotic Values for Small $n$}

Here we present the numerical values for the number of normal and tree-child networks with one,
two and three reticulation vertices. They are compared to the first and second order asymptotics.
The error is $\Ord{1/\sqrt n}$ in the first order asymptotics and
$\Ord{1/n}$ in the second order asymptotics. Thus convergence is slow. So, we chose a quadratic
scale to better visualize the convergence.

The data indicate that the more reticulation vertices the networks have, the bigger is the
constant factor in the third order term. In particular, for normal networks even the second order
asymptotics is still fairly inaccurate when the size is around 1000 vertices.

\begin{table}[ht]
	\centering
	\begin{tabular}{c||c|c|c}
		$n$ & $N_{1,n} $ & first order asymptotics  & second order asymptotics   \\
		\hline
		$7^2$
		 & $ 1.509083862\times 10^{70} $
		 & $ 2.845078723\times 10^{70}$
		 & $1.316888413 \times 10^{70} $  \\
		$9^2$
		& $ 1.424572126 \times 10^{133} $
		& $2.286221720\times 10^{133}$
	    & $1.331103718 \times 10^{133} $  \\
		$11^2$
		& $ 2.805663893\times 10^{219} $
		& $ 4.092442789\times 10^{219} $
	    & $ 2.693592858\times 10^{219} $ \\
		$13^2$
		&$  3.126424192\times 10^{330} $
		& $ 4.280475255\times 10^{330} $
		& $ 3.042449065\times 10^{330} $ \\
		$15^2$
		&$ 2.988746000\times 10^{467}  $
		& $ 3.911561797\times 10^{467} $
	    & $2.931078655\times 10^{467} $ \\
		$17^2$
		&$ 2.485340363\times 10^{631} $
		& $ 3.144767382\times 10^{631} $
		& $ 2.449229483\times 10^{631} $ \\
		$19^2$
		&$1.354821659 \times 10^{823} $
		& $ 1.669930393\times 10^{823} $
		& $ 1.339465018\times 10^{823} $ \\
		$21^2$
		&$2.903179416\times 10^{1043} $
		& $ 3.504201765\times 10^{1043}$
		& $ 2.876792390 \times 10^{1043} $ \\
		$23^2$
		&$1.222842196\times 10^{1293}$
		& $1.450922268\times 10^{1293}$
		& $ 1.213731650\times 10^{1293} $ \\
		$25^2$
		&$ 4.366393995\times 10^{1572}$
		& $ 5.107382228\times 10^{1572}$
		& $ 4.339243703\times 10^{1572} $ \\
		$27^2$
		&$5.040854939\times 10^{1882}$
		& $5.825548735\times 10^{1882}$
		& $ 5.014299556\times 10^{1882} $ \\
		$29^2$
		&$ 6.468853840\times 10^{2223}$
		& $ 7.398904501\times 10^{2223}$
		& $ 6.439612952\times 10^{2223} $ \\
		$31^2$
		& $ 2.903035924 \times 10^{2596}$
		& $ 3.290787336\times 10^{2596}$
		& $ 2.891652790\times 10^{2596} $ \\
	\end{tabular}

\vspace*{0.3cm}

	\begin{tabular}{c||c|c|c}
		$n$ & $N_{2,n} $ & first order asymptotics  & second order asymptotics \\
		\hline
		$7^2$
		& $ 1.974631541 \times 10^{72} $
		&$ 8.538792514\times 10^{72}$
		&$ -0.634169808\times 10^{72}$  \\
		$9^2$
		& $ 0.652084068 \times 10^{136} $
		& $1.874987588\times 10^{136}$
		& $ 0.308355286 \times 10^{136} $  \\
		$11^2$
		& $3.279694748\times 10^{222} $
		& $ 7.489681863\times 10^{222} $
		& $ 2.369541406\times 10^{222} $ \\
		$13^2$
		&$  0.775961070\times 10^{334} $
		& $ 1.528183172\times 10^{334} $
		& $ 0.644201521\times 10^{334} $ \\
		$15^2$
		&$ 1.393399104\times 10^{471}  $
		& $ 2.475285199\times 10^{471} $
	    & $1.234361225\times 10^{471} $ \\
		$17^2$
		&$1.993994409\times 10^{635} $
		& $ 3.283176454\times 10^{635} $
		& $ 1.830875936\times 10^{635} $ \\
		$19^2$
		&$1.751242120 \times 10^{827} $
		& $ 2.720337482\times 10^{827} $
		& $ 1.643673031\times 10^{827} $ \\
		$21^2$
		&$ 5.742638603\times 10^{1047} $
		& $ 8.518758291\times 10^{1047}$
		& $ 5.468278229 \times 10^{1047} $ \\
		$23^2$
		&$3.551625232\times 10^{1297}$
		& $ 5.075344229\times 10^{1297}$
		& $3.415952745\times 10^{1297} $ \\
		$25^2$
		&$ 1.799912695\times 10^{1577}$
		& $ 2.493838978\times 10^{1577}$
		& $1.743703701\times 10^{1577} $ \\
		$27^2$
		&$ 2.866620516\times 10^{1887}$
		& $ 3.869919308\times 10^{1887}$
		& $ 2.792091620\times 10^{1887} $ \\
		$29^2$
		&$ 4.954032473\times 10^{2228}$
		& $ 6.541380718\times 10^{2228}$
		& $4.845159005\times 10^{2228} $ \\
		$31^2$
		& $ 2.932551027\times 10^{2601}$
		& $ 3.798889014\times 10^{2601}$
		& $ 2.877366178\times 10^{2601} $ \\
	\end{tabular}

\vspace*{0.3cm}

	\begin{tabular}{c||c|c|c}
		$n$ & $N_{3,n} $ & first order asymptotics  & second order asymptotics \\
		\hline
		$7^2$
		 & $ 1.365816004\times 10^{74} $
		 &$ 17.08470069\times 10^{74}$
		 &$ -0.104456524 \times 10^{74} $  \\
		$9^2$
		& $ 1.755204956 \times 10^{138} $
		& $10.25149464\times 10^{138}$
		& $ -2.596848522\times 10^{138} $  \\
		$11^2$
		 & $ 2.360997970\times 10^{225}$
		 & $ 9.138036014\times 10^{225} $
		 & $ -0.232461034\times 10^{225} $ \\
		$13^2$
		 &$ 1.215402285 \times 10^{337} $
		 & $ 3.637203298\times 10^{337} $
		 & $ 0.481278310\times 10^{337} $ \\
		$15^2$
		 &$ 4.159928205\times 10^{474}  $
		 & $ 10.44260944\times 10^{474} $
		 & $2.589887410\times 10^{474} $ \\
		$17^2$
		&$1.034136900\times 10^{639} $
		 & $ 2.285118173\times 10^{639} $
		  & $0.768898277\times 10^{639} $ \\
		$19^2$
		&$ 1.472741571 \times 10^{831} $
		& $2.954309176\times 10^{831} $
		 & $ 1.200409327\times 10^{831} $ \\
		$21^2$
		&$0.742450513\times 10^{1052} $
		& $ 1.380613859\times 10^{1052}$
		 & $ 0.639038343\times 10^{1052} $ \\
		$23^2$
		&$6.765254066\times 10^{1301}$
		& $ 11.83574504\times 10^{1301}$
		 & $6.031172891\times 10^{1301} $ \\
		$25^2$
		 &$ 4.878745045\times 10^{1581}$
		 & $ 8.117965422\times 10^{1581}$
		 & $ 4.455195521\times 10^{1581} $ \\
		$27^2$
		&$1.074095703\times 10^{1892}$
		& $ 1.713861489\times 10^{1892}$
		& $ 0.997859210\times 10^{1892} $ \\
		$29^2$
		&$ 2.503773287\times 10^{2233}$
		& $ 3.855495246\times 10^{2233}$
		& $ 2.355863510\times 10^{2233} $ \\
		$31^2$
		& $1.957523560\times 10^{2606}$
		& $ 2.923628151\times 10^{2606}$
		& $1.859821038\times 10^{2606} $ \\
	\end{tabular}
\end{table}

\vspace*{0.3cm}

\begin{table}
    \centering
	\begin{tabular}{c||c|c|c}
		$n$ & $T_{1,n} $ & first order asymptotics  & second order asymptotics   \\
		\hline
		$7^2$
		 & $ 2.295774923\times 10^{70} $
		 & $ 2.845078723\times 10^{70}$
		 & $2.335681951 \times 10^{70} $  \\
		$9^2$
		& $ 1.948607480 \times 10^{133} $
		& $2.286221720\times 10^{133}$
	    & $1.967849052 \times 10^{133} $  \\
		$11^2$
		& $ 3.603212411\times 10^{219} $
		& $ 4.092442789\times 10^{219} $
	    & $ 3.626159479\times 10^{219} $ \\
		$13^2$
		&$  3.850668993\times 10^{330} $
		& $ 4.280475255\times 10^{330} $
		& $ 3.867799857\times 10^{330} $ \\
		$15^2$
		&$ 3.573001570\times 10^{467}  $
		& $ 3.911561797\times 10^{467} $
	    & $3.584734083\times 10^{467} $ \\
		$17^2$
		&$ 2.905589575\times 10^{631} $
		& $ 3.144767382\times 10^{631} $
		& $ 2.912921415\times 10^{631} $ \\
		$19^2$
		&$1.556662281 \times 10^{823} $
		& $ 1.669930393\times 10^{823} $
		& $ 1.559775267\times 10^{823} $ \\
		$21^2$
		&$3.289723172\times 10^{1043} $
		& $ 3.504201765\times 10^{1043}$
		& $ 3.295065305 \times 10^{1043} $ \\
		$23^2$
		&$1.370016200\times 10^{1293}$
		& $1.450922268\times 10^{1293}$
		& $ 1.371858728\times 10^{1293} $ \\
		$25^2$
		&$ 4.845849891\times 10^{1572}$
		& $ 5.107382228\times 10^{1572}$
		& $ 4.851336052\times 10^{1572} $ \\
		$27^2$
		&$5.549770266\times 10^{1882}$
		& $5.825548735\times 10^{1882}$
		& $ 5.555132339\times 10^{1882} $ \\
		$29^2$
		&$ 7.073239889\times 10^{2223}$
		& $ 7.398904501\times 10^{2223}$
		& $ 7.079140649\times 10^{2223} $ \\
		$31^2$
		& $ 3.155446557 \times 10^{2596}$
		& $ 3.290787336\times 10^{2596}$
		& $ 3.157742486\times 10^{2596} $ \\
	\end{tabular}

\vspace*{0.3cm}

	\begin{tabular}{c||c|c|c}
		$n$ & $T_{2,n} $ & first order asymptotics  & second order asymptotics  \\
		\hline
		$7^2$
		& $ 4.640516422 \times 10^{72} $
		&$ 8.538792514\times 10^{72}$
		&$ 5.481138402\times 10^{72}$  \\
		$9^2$
		& $ 1.224486692 \times 10^{136} $
		& $1.874987588\times 10^{136}$
		& $ 1.352776820 \times 10^{136} $  \\
		$11^2$
		& $5.411620540\times 10^{222} $
		& $ 7.489681863\times 10^{222} $
		& $ 5.782968374\times 10^{222} $ \\
		$13^2$
		&$  1.176422691\times 10^{334} $
		& $ 1.528183172\times 10^{334} $
		& $ 1.233522622\times 10^{334} $ \\
		$15^2$
		&$ 1.989643302\times 10^{471}  $
		& $ 2.475285199\times 10^{471} $
	    & $2.061643874\times 10^{471} $ \\
		$17^2$
		&$ 2.722734486\times 10^{635} $
		& $ 3.283176454\times 10^{635} $
		& $ 2.799076281\times 10^{635} $ \\
		$19^2$
		&$2.309772743 \times 10^{827} $
		& $ 2.720337482\times 10^{827} $
		& $ 2.361449330\times 10^{827} $ \\
		$21^2$
		&$7.367315911\times 10^{1047} $
		& $ 8.518758291\times 10^{1047}$
		& $ 7.501931598 \times 10^{1047} $ \\
		$23^2$
		&$4.454475775\times 10^{1297}$
		& $ 5.075344229\times 10^{1297}$
		& $4.522213731\times 10^{1297} $ \\
		$25^2$
		&$ 2.215316571\times 10^{1577}$
		& $ 2.493838978\times 10^{1577}$
		& $2.243793885\times 10^{1577} $ \\
		$27^2$
		&$ 3.472411188\times 10^{1887}$
		& $ 3.869919308\times 10^{1887}$
		& $ 3.510643410\times 10^{1887} $ \\
		$29^2$
		&$ 5.919519663\times 10^{2228}$
		& $ 6.541380718\times 10^{2228}$
		& $5.975973477\times 10^{2228} $ \\
		$31^2$
		& $ 3.462830748\times 10^{2601}$
		& $ 3.798889014\times 10^{2601}$
		& $ 3.491714733\times 10^{2601} $ \\
	\end{tabular}

\vspace*{0.3cm}

	\begin{tabular}{c||c|c|c}
		$n$ & $T_{3,n} $ & first order asymptotics  & second order asymptotics   \\
		\hline
		$7^2$
		 & $ 4.905522940\times 10^{74} $
		 &$ 17.08470069\times 10^{74}$
		 &$ 7.907916294 \times 10^{74} $  \\
		$9^2$
		& $ 4.488502332 \times 10^{138} $
		& $10.25149464\times 10^{138}$
		& $ 5.968713581\times 10^{138} $  \\
		$11^2$
		 & $ 4.976717574\times 10^{225}$
		 & $ 9.138036014\times 10^{225} $
		 & $ 6.014536995\times 10^{225} $ \\
		$13^2$
		 &$ 2.258644701\times 10^{337} $
		 & $ 3.637203298\times 10^{337} $
		 & $ 2.585228301\times 10^{337} $ \\
		$15^2$
		 &$ 7.072035184\times 10^{474}  $
		 & $ 10.44260944\times 10^{474} $
		 & $7.825035424\times 10^{474} $ \\
		$17^2$
		&$ 1.645107613\times 10^{639} $
		 & $ 2.285118173\times 10^{639} $
		  & $1.779711540\times 10^{639} $ \\
		$19^2$
		&$ 2.225232465 \times 10^{831} $
		& $2.954309176\times 10^{831} $
		 & $ 2.369675891\times 10^{831} $ \\
		$21^2$
		&$1.076588119\times 10^{1052} $
		& $ 1.380613859\times 10^{1052}$
		 & $ 1.133422020\times 10^{1052} $ \\
		$23^2$
		&$9.485462012\times 10^{1301}$
		& $ 11.83574504\times 10^{1301}$
		 & $9.900887645\times 10^{1301} $ \\
		$25^2$
		 &$ 6.651391940\times 10^{1581}$
		 & $ 8.117965422\times 10^{1581}$
		 & $ 6.897042119\times 10^{1581} $ \\
		$27^2$
		&$1.430044284\times 10^{1892}$
		& $ 1.713861489\times 10^{1892}$
		& $ 1.475194062\times 10^{1892} $ \\
		$29^2$
		&$ 3.266427497\times 10^{2233}$
		& $ 3.855495246\times 10^{2233}$
		& $ 3.355617999\times 10^{2233} $ \\
		$31^2$
		& $2.509177651\times 10^{2606}$
		& $ 2.923628151\times 10^{2606}$
		& $2.569025778\times 10^{2606} $ \\
	\end{tabular}
\end{table}

\end{document}